\documentclass[preprint,12pt]{elsarticle}

\usepackage{amssymb}
\usepackage{bm}
\usepackage{amsmath}
\journal{Journal of Multivariate Analysis
}

\newtheorem{thm}{Theorem}
\newtheorem{prop}{Proposition}
\newtheorem{lem}{Lemma}
\newdefinition{rmk}{Remark}
\newproof{pf}{Proof}
\newproof{pot}{Proof of Theorem \ref{thm:thm31}}
\newproof{pop1}{Proof of Proposition \ref{prop:prop31}}
\newproof{pop2}{Proof of Proposition \ref{prop:prop32}}
\newproof{pop3}{Proof of Proposition \ref{prop:prop33}}
\newproof{pol1}{Proof of Lemma \ref{lem:lem71}}
\newproof{pol2}{Proof of Lemma \ref{lem:lem72}}
\newproof{pol3}{Proof of Lemma \ref{lem:lem73}}
\newcommand{\RP}{{\rm P}}

\newcommand{\RE}{{\rm E}}

\newcommand{\sign}{{\rm sign}}
\newcommand{\argmin}{\mathop{\rm argmin}\limits}
\newcommand{\bxi}{\bm{X}_i(t)}
\newcommand{\bwi}{\bm{W}_i(t)}
\newcommand{\bgt}{\bm{g}(t)}
\newcommand{\obt}{\overline{\bm{B}}(t)}
\newcommand{\bgmj}{\bm{\gamma}_{-1j}}
\newcommand{\bgj}{\bm{\gamma}_j}
\newcommand{\bgg}{\bm{\gamma}}
\newcommand{\bt}{\bm{\theta}}
\newcommand{\wbt}{\widehat{\bm{\theta}}}
\newcommand{\wt}{\widehat{\theta}}
\newcommand{\wbb}{\widehat{\bm{b}}}
\begin{document}

\begin{frontmatter}


\title{Variable selection and structure identification for varying coefficient
Cox models}
%\tnoteref{label1}
\author{Toshio Honda\corref{cor1}}
\ead{t.honda@r.hit-u.ac.jp}
\address{Graduate School of Economics, Hitotsubashi University,
Kunitachi, Tokyo 186-8601, Japan}
\cortext[cor1]{Corresponding author}
%%%%%%%%%%%%%%%%%%%%%%%%% Prof. Yabe
\author{Ryota Yabe\corref{cor2}
%\fnref{fn1}
}
\ead{ryotayabe@shinshu-u.ac.jp}
\address{Department of Economics, Shinshu University,
Matsumoto, Nagano 390-8621, Japan}


%%%%%%%%%%%%%%%%%%%%%%%%%%%%%%%%%%%%%%%%%%%%%%%%%%%%%%%%%%%%%%%%% Abstract
\begin{abstract}
We consider varying coefficient Cox models with high-dimensional covariates.
We apply the group Lasso method to these models and propose a variable selection
procedure. Our procedure copes with variable selection and structure identification
from a high dimensional varying coefficient model to a semivarying coefficient model
simultaneously. We derive an oracle inequality and closely examine restrictive eigenvalue
conditions, too. In this paper, we give the details for Cox models with time-varying
coefficients. The theoretical results on variable selection can be easily extended to
some other important models and we briefly mention those models since those models can be
treated in the same way. The models considered in this paper are the most popular
models among structured nonparametric regression models.
The results of a small numerical study are also given.
\end{abstract}


\end{frontmatter}

% \linenumbers
%%%%%%%%%%%%%%%%%%%%%%%%%%%%%%%%%%%%%%%%%%%%%%%%%%%%%%%%%%%%% Main text
\section{Introduction}
\label{intro}
%\begin{enumerate}
%\item our model and references on varying coefficient Cox models, briefly our problem
%\item motivation for high-dimensinal problems and general references on the problems
%\item references on high-dimensional Cox models
%\item our problem and related references (SCAD)
%\end{enumerate}

The Cox model is one of the most popular and useful models to analyze censored survival
data. Since the Cox model was proposed in Cox\cite{Cox1972},
many authors have studied a lot of extensions 
or variants of the original Cox model to deal with complicated situations or carry out
more flexible statistical analysis. In this paper, we consider varying coefficient models
and additive models with high-dimensional covariates. These models with moderate numbers
of covariates are investigated in many papers, for example, Huang et al.\cite{HKST2000},
Cai and Sun\cite{CS2003}, and Cai et al.\cite{CFZZ2007}.

We apply the group Lasso (for example, see Lounici et al.\cite{LPGT2011}) to
varying coefficient models with high-dimensional covariates to carry out variable
selection and structure identification simultaneously. Although
we focus on time-varying coefficient
models here, our method can be applied to variable selection for another type of varying
coefficient models and additive models and we briefly mention how to apply
and how to derive the theoretical results.

Suppose that we observe censored survival times $T_i$ and high-dimensional
covariates $\bxi=(X_{i1}(t), \ldots, X_{ip}(t))^T$. More specifically, we have
$n$ i.i.d. observations of
\begin{equation}
T_i = \min \{ T_{0i}, C_i \}, \qquad \delta_i =I\{ T_{0i} \le C_i\},
\label{eqn:e101}
\end{equation}
and $p$-dimensional covariate $\bxi$ on the time interval $[0,\tau]$, where
$T_{0i}$ is an uncensored survival time and $C_i$ is a censoring time
satisfying the condition of the independent censoring mechanism as in
section 6.2 of Kalbfleisch and Prentice\cite{KP2002}. Hereafter we set
$\tau =1$ for simplicity of presentation. Note that $p$ can be very
large compared to $n$ in this paper, for example, $ p = O(n^{c_p})$ for
a very large positive constant $c_p$ or $p=O( \exp( n^{c_p} ) )$
for a sufficiently small positive constant $c_p$.
We assume that the standard setup
for the Cox model holds as in chapter 5 of \cite{KP2002} and that $T_i$ or
$N_i(t)=I\{ t \ge T_i\}$ has the following compensator $\Lambda_i(t)$
with respect to a suitable filtration $\{ {\cal F}_t\}$:
\begin{equation}
d\Lambda_i(t) = Y_i(t) \exp \{ \bxi^T\bgt \} \lambda_0 (t) dt,
\label{eqn:e103}
\end{equation}
where $Y_i (t) = I \{ t \le T_i \}$, $\bgt = (g_1(t),
\ldots, g_p(t))^T$ is a vector of unknown functions on $[0,1]$,
$\bm{a}^T$ denotes the transpose of $\bm{a}$, and
$\lambda_0 (t) $ is a baseline hazard function. As in chapter 5
of \cite{KP2002}, $\bxi$ is predictable and
\begin{equation}
M_i(t) = N_i(t) - \Lambda_i (t) 
\label{eqn:e105}
\end{equation}
is a martingale with respect to $\{ {\cal F}_t\}$. In the
original Cox model, $\bgt$ is a vector of constants and we estimate this
constant coefficient vector by maximizing the partial likelihood.

% high-dimensional problems
In this paper, we are interested in estimating $\bgt$ in (\ref{eqn:e103}).
Recently we have many cases where there are (ultra) high-dimensional
covariates due to drastic development of data collecting technology.
In such high-dimensional data,
usually only a small part of covariates are relevant.
However, we cannot directly apply standard or traditional estimating procedures
to such high-dimensional data. Thus now a lot of methods for variable selection are available, for example,
SCAD and Lasso procedures. See B\"uhlmann and van de Geer\cite{BG2011} and
Hastie et al.\cite{HTW2015} for excellent reviews of these procedures for variable
selection. See also Bickel et al.\cite{BRT2009} and Zou\cite{Zou2006} for the Lasso
and the adaptive Lasso, respectively.

% high-dimensional original Cox model
As for high dimensional Cox models with constant coefficient, 
Bradic et al.\cite{BFJ2011} studied the SCAD method and
Huang et al.\cite{HSYYZ2013} considered the Lasso procedure.
The authors of \cite{HSYYZ2013} developed new ingenious
techniques to derive oracle inequalities. We will fully use their
techniques to derive our theoretical results such as an oracle inequality.
In addition, Zhang and Luo\cite{ZL2007} proposed an adaptive Lasso estimator
for the Cox model. Some variable screening procedures have also been proposed in
Zhao and Li\cite{ZL2012} and Yang et al.\cite{YYLB2016}, to name just a few.

% our purpose and closely related papers
In this paper, we propose a group Lasso procedure
to select relevant covariates and identify the covariates with constant
coefficients among the relevant covariates, namely the
true semivarying coefficient model from the original varying coefficient
model. We can achieve this goal by the proposed group Lasso with a suitable
threshold value or a two-stage procedure consisting of the proposed
one and an adaptive Lasso procedure as in Yan and Huang\cite{YH2012}
and Honda and H\"ardle\cite{HH2014}. In \cite{YH2012},
the authors proposed an adaptive Lasso procedure for structure identification
with no theoretical result. Our procedure can be applied to
the varying coefficient model with an index variable $Z_i(t)$:
\begin{equation}
d\Lambda_i(t) = Y_i(t) \exp \{ g_0(Z_i(t)) + \bxi^T\bm{g}(Z_i(t)) \}
\lambda_0 (t) dt
\label{eqn:e107}
\end{equation}
and the additive model:
\begin{equation}
d\Lambda_i(t) = Y_i(t) \exp \Big\{ \sum_{j=1}^p
g_j(X_{ij}(t)) \Big\} \lambda_0 (t) dt.
\label{eqn:e109}
\end{equation}
We mention these model later in section \ref{sec:other}.

Some authors considered the same problem by using SCAD.
For example, see Lian et al.\cite{LLL2013} and Zhang et al.\cite{ZWL2014}.
They proved the existence of local optimizer satisfying the same convergence
rate as ours. In contrast, we prove the existence of the global solution
with desirable properties. In Bradic and Song\cite{BS2015}, the authors
applied penalties similar to ours to additive models and obtained
theoretical results in another complicated manner. 
We have derived a better convergence rate for our procedures to varying coefficient models
by exploiting the martingale structure very carefully under much simpler assumptions
given in section \ref{sec:group}. See Remark 1 in section \ref{sec:oracle} about
the convergence rate.
% Besides, their Lemma 3 seems to fail to hold for our procedures and basis functions.
We also carefully examined the
RE (restrictive eigenvalue) conditions. While the other authors considered
the $L_2$ norm of the estimated second derivatives for additive models,
we adopt the orthogonal decomposition approach. We give some details on why we
have adopted the orthogonal decomposition approach in \ref{sec:bspline}.
% rate improvement by dealing specific cases, a theoretical
% problem in Lemma 3, and careful examination of RE conditions.

%%%% organization of this paper
% section variable : description of our procedure
% section oracle : oracle inequality and RE conditions
% section other : other two models
% section numerical : if any
% section proofs :
% section appendix : supplementary?  

This paper is organized as follows. In section \ref{sec:group},
we describe our group Lasso procedure for time-varying coefficient
models. Then we present our theoretical results in section \ref{sec:oracle}.
We mention the two other models in section \ref{sec:other}. The results
of a small simulation study are given in section \ref{sec:numerical}.
The proofs of our theoretical results are postponed to section \ref{sec:proofs}
and section \ref{sec:conclusion} concludes this paper.
We collected useful properties of our basis functions and the proofs of technical
lemmas in Appendices A-E.

We define some notation and symbols here.
In this paper, $C$, $C_1$, $C_2$, $\ldots$ are positive generic constants
and their values change from line to line. For a vector $\bm{a}$,
$|\bm{a}|$, $|\bm{a}|_1$,
and $|\bm{a}|_\infty$ mean the $L_2$ norm, the $L_1$ norm, and the sup norm,
respectively. For a function $g$ on $[0,1]$, $\| g \|$, $\| g \|_1$, and $\|
g \|_\infty$ stand for the $L_2$ norm, the $L_1$ norm, and the sup norm,
respectively. For a symmetric matrix $A$, we denote the
minimum and maximum eigenvalues by $\lambda_{\min}(A)$
and $\lambda_{\max}(A)$, respectively. Besides, $\sign (a)$ is the
sign of a real number $a$ and $a_n \sim b_n$ means
there are positive constants $C_1$ and $C_2$ such that
$ C_1 < a_n/b_n < C_2 $.
We write $\overline{\cal S}$ for the complement of a set ${\cal S}$.

%\begin{equation}\label{eqn:e101}\end{equation}

\section{Group Lasso procedure}
\label{sec:group}
First we decompose $g_j(t)$, $j=1,\ldots, p$, into the constant part
and the non-constant part:
\begin{equation}
g_j(t) = g_{cj} + g_{nj}(t),
\label{eqn:e201}
\end{equation}
where $\int_0^1 g_{nj}(t) dt =0$. When $g_j(t) \not \equiv 0$, $g_j(t)$
is a non-zero constant or a non-constant function. We denote the index sets
of relevant covariates by
\begin{equation}
{\cal S}_c = \{ j\, | \, g_{cj}\neq 0 \}\quad
{\rm and} \quad {\cal S}_n = \{ j\, | \, g_{nj}(t) \not \equiv 0 \}
\label{eqn:e203}
\end{equation}
and set
\[
s_c = \# {\cal S}_c, \quad  s_n = \# {\cal S}_n,\quad
{\rm and}\quad s_o= s_c + s_n,
\]
where $\# A$ is the number of the elements of a set $A$. We implicitly
assume that $s_0$ is bounded or much smaller than $n$. Besides, we assume
\begin{equation}
{\cal S}_n \subset {\cal S}_c.
\label{eqn:e205}
\end{equation}
We may incidentally have $g_{cj}=0$ for $j\in {\cal S}_n$. However,
this will rarely happen and $g_{cj}$ should be free if $ g_{nj}(t)
\not \equiv 0 $.

Next we introduce our spline basis $\overline{\bm{B}}(t)$ to
approximate $g_j(t)$, $j=1,\ldots, p$. We construct
$\obt$ from the $L$-dimensional equispaced B-spline basis
$\bm{B}_0(t)=(b_{01}(t), \ldots, b_{0L}(t))^T$ on $[0,1]$ 
and the basis has the following properties :
\begin{equation}
\obt =
\begin{pmatrix}
b_1(t)\\
b_2(t)\\
\vdots\\
b_L(t)
\end{pmatrix}
=
\begin{pmatrix}
1/\sqrt{L}\\
\bm{B}(t)
\end{pmatrix}
= A_0 \bm{B}_0(t) 
\quad {\rm and}\quad
\int_0^1 \obt \overline{\bm{B}}^T(t) dt= L^{-1}I,
\label{eqn:e207}
\end{equation}
where
\[
A_0=
\begin{pmatrix}
\bm{a}_{01}^T \\
\bm{a}_{02}^T \\
\vdots \\
\bm{a}_{0L}^T \\
\end{pmatrix}
=
\begin{pmatrix}
\bm{1}^T/\sqrt{L}\\
A_{-1}
\end{pmatrix}\quad {\rm (say)}
\]
and $ \bm{1}=(1,\ldots, 1)^T$.
Note that for $j=1, \ldots, L$,
\[
b_j(t) = \bm{a}_{0j}^T  \bm{B}_0(t)
\]
and that $1/\sqrt{L}$ and $\bm{B}(T)=(b_2(t), \ldots , b_L(t))^T$
in (\ref{eqn:e207}) are designed for $g_{cj}$ and $g_{nj}(t)$, respectively.
Recall that $\bm{1}^T\bm{B}_0(t) \equiv 1$ and see
Schumaker\cite{Schumaker2007} for the definition of B-spline bases.
We have collected how to construct $\obt$ and $A_0$ and some useful properties of
$\obt$ and $A_0$ in \ref{sec:basis}. We can use another basis which has
desirable properties such as (\ref{eqn:e901}), (\ref{eqn:e906}),
and (\ref{eqn:e904}) in \ref{sec:basis}.

We impose some technical assumptions on $\bm{g}(t)$.

%%%%%%%%%%%%%%%% Assumption G
\noindent
{\bf Assumption G :} $g_j(t)$, $j=1,\ldots, p$, are twice continuously
differentiable and there is a positive constant $C_g$ such that
\[
\sum_{j=1}^p\| g_j \|_\infty \le C_g, \quad
\sum_{j=1}^p\| g_j' \|_\infty \le C_g, \quad {\rm and}\quad
\sum_{j=1}^p\| g_j'' \|_\infty \le C_g.
\]
Besides we have
\[
\min_{j \in {\cal S}_c} |g_{cj}|L^2 \to \infty\quad
{\rm and}\quad \min_{j \in {\cal S}_n} \|g_{nj}\| L^2 \to \infty .
\]

\medskip
Hereafter we take $L=c_L n^{1/5}(c_L>0)$ for simplicity of presentation
and the order of the B-spline basis should be larger than or equal to 2.
The latter of Assumption G means relevant coefficient functions are
larger than the approximation error. As for the identifiability of
$\bm{g}(t)$, we need an assumption such as $\lambda_{\min}( \RE\{
\overline \Sigma\} ) > C_1/L$ for a positive constant $C_1$,
where $ \RE \{ \overline \Sigma \} $ is defined in Proposition \ref{prop:prop33}.

When Assumption G holds, there are $\bm{\gamma}_j^*
=(\gamma_{1j}^*, \bm{\gamma}_{-1j}^{*T})^T\in R^L $, $j=1,\ldots, p$,
such that for a positive constant $C_{approx}$ depending on $C_g$,
\begin{equation}
\sum_{j=1}^p\| g_j - \obt^T \bgj^* \|_\infty \le C_{approx}L^{-2}.
\label{eqn:e208}
\end{equation}
When $j \in {\cal S}_c$, we can take $\gamma_{1j}^* = \sqrt{L}g_{cj}$
and $\bm{\gamma}_{-1j}^{*}\in R^{L-1}$ depends on $g_{jn}(t)$.
If $j \in \overline{\cal S}_n$, we take $\bm{\gamma}_{-1j}^{*}=0$.
When $j \in \overline{\cal S}_c$, we set $\bgj^* =0$. See \ref{sec:basis}
for more details on these $\bm{\gamma}_j^*
=(\gamma_{1j}^*, \bm{\gamma}_{-1j}^{*T})^T$.

We state assumptions on our Cox model before we describe the log partial
likelihood for new covariates
\begin{equation}
\bm{W}_i(t) = \bxi \otimes \obt ,
\label{eqn:e209}
\end{equation}
where $\otimes$ means the Kronecker product.

%%%%%%%%%%%%%%%% Assumption M
\noindent
{\bf Assumption M :} $|X_{1j}(t)|\le C_X$ uniformly in $j$ and $t$
for a positive constant $C_X$. We also have $\RE \{ Y_1(1) \}\ge C_Y$
for a positive constant $C_Y$. Besides, the baseline hazard
function is bounded from above and satisfies
$\lambda_0(t) \ge C_\lambda$ on $[0,1]$ for a positive
constant $C_\lambda$.

\medskip
The first one is used to evaluate the inside of the exponential function
and the other ones are standard in the literature.

%%%%%%%%%%%%%%%%%%%%%%%%% Partial Likelihood
We denote the log partial likelihood by $L_p(\bm{\gamma})$ :
\begin{equation}
L_p(\bm{\gamma}) =\frac{1}{n}\sum_{i=1}^n \int_0^1
\bm{\gamma}^T\bwi d N_i(t) - \int_0^1
\log \Big[ \sum_{i=1}^n Y_i(t) \exp \{ \bm{\gamma}^T\bwi \}
\Big] d \overline N (t),
\label{eqn:e211}
\end{equation}
where
$\bm{\gamma} = ( \bm{\gamma}_1^T, \ldots, \bm{\gamma}_p^T)^T\in 
R^{pL}$ and
$\overline N (t) = n^{-1}\sum_{i=1}^n N_i(t)$. We also
use the same sample mean notation for $M_i(t)$ and $Y_i(t)$.

Set
\begin{equation}
l_p( \bm{\gamma} ) = -L_p ( \bm{\gamma} )
\label{eqn:e212}
\end{equation}
for notational convenience. Then we should minimize this $l_p( \bm{\gamma} ) $
with respect to $ \bm{\gamma} $. However, when $pL$ is larger than $n$,
we cannot carry out this minimization properly and we add some penalty as
in the literature on high-dimensional data. We define two convex penalties here :
\begin{align}
P_1(\bm{\gamma}) & = \sum_{j=1}^p ( |\gamma_{1j}| + | \bgmj |)
\label{eqn:e213}\\
\intertext{and}
P_h(\bm{\gamma}) & = \sum_{j=1}^p ( |\gamma_{1j}|^q + | \bgmj |^q)^{1/q}
+  \sum_{j=1}^p | \bgmj | 
\label{eqn:e215}
\end{align}
for some $q>1$.

This $P_1(\bm{\gamma})$ plays the role of the $L_1$ norm for $\bm{\gamma}
\in R^{pL}$ and is a very important technical tool in this paper. Besides,
we define a kind of sup norm $ P_\infty (\bm{\gamma}) $ by
\begin{equation}
P_\infty (\bm{\gamma}) = \max_{1\le j \le p}|\gamma_{1j}| \vee |\bgmj|,
\label{eqn:e217}
\end{equation}
where $a\vee b = \max \{ a, b \} $. This is also an important tool.

We defined the penalty in (\ref{eqn:e215}) by taking the assumption
in (\ref{eqn:e205}) into consideration and following
Zhao et al.\cite{ZRY2009} and Zhao and Leng\cite{ZL2016}. Thus
our group Lasso objective functions are
\begin{equation}
Q_h(\bm{\gamma}; \lambda ) =l_p(\bm{\gamma}) + \lambda P_h (\bm{\gamma})
\quad {\rm and}\quad
Q_1(\bm{\gamma}; \lambda ) =l_p(\bm{\gamma}) + \lambda P_1 (\bm{\gamma}).
\label{eqn:e219}
\end{equation}

Our group Lasso estimate is given by
\[
\widehat \bgg = \argmin_{\bgg \in R^{pL}}Q_h (\bgg ; \lambda)
\quad {\rm or}\quad
\widehat \bgg = \argmin_{\bgg \in R^{pL}}Q_1 (\bgg ; \lambda).
\]

If we are interested in only variable selection, we should
minimize
\begin{equation}
Q(\bm{\gamma}; \lambda )=l_p(\bm{\gamma}) +
\lambda \sum_{j=1}^p|\bm{\gamma}_j|.
\label{eqn:e231}
\end{equation}

The KKT condition implies that for $a=h$ or $1$,
\begin{equation}
\frac{\partial l_p}{\partial \bgg_j}
(\widehat \bgg) = -\lambda \nabla_j
P_a ( \widehat \bgg ), \quad j=1, \ldots, p,
\label{eqn:e220}
\end{equation}
where $ \nabla_j P_a (\bgg) $ is the subgradient of $P_a (\bgg)$
with respect to $\bgj$. See chapter 5 of \cite{HTW2015} about convex
optimality conditions. We give explicit expressions of these subgradients
in \ref{sec:subgradients} for reference.
Consequently from (\ref{eqn:e207}), our estimates
of $g_{cj}$ and $g_{nj}$ are
\begin{equation}
\widehat g_{cj} = \widehat \gamma_{1j}/\sqrt{L}
\quad {\rm and} \quad \widehat g_{nj}(t) = \bm{B}^T(t) \widehat \bgg_{-1j}.
\label{eqn:e218}
\end{equation}

If we choose a threshold value $t_\lambda$ based on our theoretical results
in section \ref{sec:oracle} and define $\widehat {\cal S}_c$ and
$\widehat {\cal S}_n$ by
\begin{equation}
\widehat {\cal S}_c= \{ j\, | \, | \widehat g_{cj} |> t_\lambda\}
\quad {\rm and}\quad
\widehat {\cal S}_n= \{ j\, | \, \| \widehat g_{nj} \| > t_\lambda\},
\label{eqn:e221}
\end{equation}
they are consistent estimators of ${\cal S}_c$ and
${\cal S}_n$, respectively. Or we can apply an adaptive
Lasso procedure to estimate the true semivarying coefficient model.

We state our theoretical results only for $Q_h(\bgg ; \lambda)$
in section \ref{sec:oracle} since we can deal with $Q_1(\bgg ; \lambda)$
and $Q(\bm{\gamma}; \lambda )$
in the same way. In terms of numerical optimization, $Q_1(\bgg ; \lambda)$
seems to be more tractable and we focused on $Q_1(\bgg ; \lambda)$
in our numerical study. When the group Lasso based on
$Q_1(\bgg ; \lambda)$ concludes that $  \| g_{nj} \| >0 $
and $ | g_{cj} |=0 $, we should take (\ref{eqn:e205}) into consideration and 
modify this conclusion to the one that both of them are relevant
for this $j$.

%\smallskip
%\medskip
%\begin{equation}\label{eqn:e2}\end{equation}

\section{Oracle inequality}
\label{sec:oracle}
An oracle inequality for $\widehat{\bm{\gamma}}$
from $Q_h(\bm{\gamma}; \lambda )$
is given in Theorem \ref{thm:thm31}.
First we define some notation. We borrow some notation from
\cite{HSYYZ2013} and proceed as in \cite{HSYYZ2013}.
Some other notation is standard in the literature of the Cox model
and the Lasso.

Let $\bgg_{\cal S}$ consist of $\{\gamma_{1j} \}_{j\in {\cal S}_c}$ and
$\{\bgmj \}_{j\in {\cal S}_n}$. On the other hand, $\bgg_{\overline{\cal S}}$
consists of $\{\gamma_{1j} \}_{j\in \overline{\cal S}_c}$ and
$\{\bgmj \}_{j\in \overline{\cal S}_n}$.  
%As in \cite{HSYYZ2013}, \begin{align}
%\begin{equation}
%\dot{l}_p (\bgg) 
%&
%=,\quad
%{\rm and}\quad 
%
%}(\bgg )
%,\quad {\rm and}\nonumber \\
%D^S(\widehat \bgg, \bgg ) &  = ( \widehat \bgg - \bgg )^T (
%\dot{l}( \widehat \bgg ) - \dot{l}( \bgg ) ) \ge 0,
%\label{eqn:e301}
%\end{equation}
%\end{align}
%where $D^S( \widehat \bgg, \bgg ) $ is called the symmetric Bregman divergence
%and its non-negativity follows from the convexity of $l_p(\bgg )$.

We need some notation to give explicit expressions of
the derivatives of $l_p(\bgg)$.
%$\dot{l}_p(\bgg )$ and $\ddot{l}_p(\bgg )$.
\[
S^{(k)}(t,\bgg ) =\frac{1}{n}\sum_{i=1}^n Y_i(t)
\bm{W}_i^{\otimes k}(t) \exp \{ \bm{W}_i^T(t) \bgg \},
\]
where $\bm{a}^{\otimes 0}=1$, $\bm{a}^{\otimes 1}=\bm{a}$,
and $\bm{a}^{\otimes 2}=\bm{a}\bm{a}^T$. In addition,
\begin{equation}
\widetilde{\bm{W}}_n (t,\bgg ) = \frac{S^{(1)}(t,\bgg )}{ S^{(0)}(t,\bgg )}
\quad {\rm and} \quad
V_n (t,\bgg ) = \frac{S^{(2)}(t,\bgg ) }{S^{(0)}(t,\bgg ) }
- ( \widetilde{ \bm{W}}_n(t,\bgg ) )^{\otimes 2}.
\label{eqn:e301}
\end{equation}

Hence we have the following expressions of the derivatives of $l_p(\bgg)$,
which are denoted by $\dot{l}_p(\bgg )$ and $\ddot{l}_p(\bgg )$ :
\begin{align}
\frac{\partial l_p}{\partial \bgg}(\bgg ) & =-\frac{1}{n}\sum_{i=1}^n
\int_0^1\{ \bwi - \widetilde{ \bm{W}}_n(t,\bgg ) \}dN_i(t)
= \dot{l}_p(\bgg )\quad  (say)
\label{eqn:e302} \\
\intertext{and}
\frac{\partial^2 l_p}{\partial \bgg \partial \bgg^T}(\bgg )
& = \int_0^1 V_n (t,\bgg ) d\overline{N}(t)= \ddot{l}_p(\bgg )
\quad (say).
\label{eqn:e303}
\end{align}

In Proposition \ref{prop:prop31}, we prove that $\widehat \bgg$
is in a restricted parameter space. We define some
more notation to state Proposition \ref{prop:prop31}.
Set
\begin{equation}
D_l= P_\infty ( \dot{l}_p(\bgg^*) )\quad
{\rm and}\quad
\widehat{\bm{\theta}} = \widehat{ \bgg} - \bgg^*.
\label{eqn:e305}
\end{equation}
We evaluate $D_l$ later in Proposition \ref{prop:prop32}.
We define $\bt_{\cal S}$ and $\bt_{\overline{\cal S}}$
in the same way as  $\bgg_{\cal S}$ and $\bgg_{\overline{\cal S}}$.
Recall that $\bgg^* = (\bgg_1^{*T}, \ldots, \bgg_p^{*T})^T$
is given in (\ref{eqn:e208}). This proposition follows from only (\ref{eqn:e220}).

%%%%%%%%%%%%%%%%%%%%%%%%%%%%%%%%%%%%%%%%%%%%%%%%%% Prop.1
\begin{prop}\label{prop:prop31} If $\lambda > D_l$, we have
\begin{align*}
(\widehat{\bgg} - \bgg^*)^T \{ \dot{l}_p (\widehat{\bgg})
- \dot{l}_p(\bgg^*) \} & \le 
(2\lambda +D_l) P_1( \widehat{\bt}_{\cal S})
- (\lambda - D_l)P_1( \widehat{\bt}_{\overline{\cal S}})\\
\intertext{and}
(\lambda - D_l)P_1( \widehat{\bt}_{\overline{\cal S}})
%& \le D^S(\widehat{\bgg}, \bgg^*) +
% (\lambda - D_l)P_1( \widehat{\bt}_{\overline{\cal S}}) \\
& \le (2\lambda +D_l) P_1( \widehat{\bt}_{\cal S}).
\end{align*}
Therefore if $D_l \le \xi \lambda\, (\xi <1 )$, we have
\[
P_1( \widehat{\bt}_{\overline{\cal S}}) \le
\frac{2+\xi}{1-\xi} P_1( \widehat{\bt}_{\cal S}).
\] 
\end{prop}

We define a restricted parameter space $\Theta (\zeta) $ by
\[
\Theta (\zeta) = \{ \bt \in R^{pL}\, | \,
P_1( \bt_{\overline{\cal S}}) \le \zeta P_1( \bt_{\cal S}) 
% \ {\rm and }\ \bt \ne 0 
\}.
\]
For $\bt \in \Theta (\zeta) $, we have
\begin{equation}
P_1(\bt ) \le (1+\zeta )P_1( \bt_{\cal S})\quad {\rm and}\quad
P_1( \bt_{\cal S}) \le s_0^{1/2}|\bt_{\cal S}|\le
s_0^{1/2}|\bt|.
\label{eqn:e307}
\end{equation}
Recall that $s_0$ is defined just after (\ref{eqn:e203}).

To state the compatibility and restrictive eigenvalue conditions, we define
$\kappa (\zeta, \Sigma )$ and $RE (\zeta, \Sigma )$
for an n.n.d.(non-negative definite) matrix $\Sigma$ with
some modifications adapted to our setup.
\[
\kappa (\zeta, \Sigma )=\inf_{\bt \in \Theta (\zeta),\, \bt \ne 0}
\frac{s_0^{1/2}(\bt^T \Sigma \bt )^{1/2}}{P_1( \bt_{\cal S})}
\quad {\rm and}\quad 
RE (\zeta, \Sigma )=\inf_{\bt \in \Theta (\zeta),\, \bt \ne 0}
\frac{(\bt^T \Sigma \bt )^{1/2}}{|\bt |}.
\]
The latter is more commonly used in the literature of the Lasso.
It is known that 
%According to Lemma 4.1 in \cite{HSYYZ2013}, we have
\[
\kappa^2 (\zeta, \Sigma ) \ge RE^2 (\zeta, \Sigma ) \ge 
\lambda_{\min}(\Sigma)
\]
and that if $\Sigma_1 - \Sigma_2 $ is n.n.d., we also have
\[
\kappa (\zeta, \Sigma_1 )\ge
\kappa (\zeta, \Sigma_2 )\quad {\rm and}\quad
RE(\zeta, \Sigma_1 )\ge
RE (\zeta, \Sigma_2 ).
\]

Some more notation is necessary for Theorem \ref{thm:thm31}.
Set
%$C_W= 2C_X ( \lambda_{\max}(A_0A_0^T) )^{1/2}$,
%\begin{equation}
\begin{align}
C_W & = 2C_X \{ \lambda_{\max}(A_0A_0^T) \}^{1/2}, \quad
RE^*= RE\Big( \frac{2+\xi}{1-\xi}, \ddot{l}_p(\bgg^*) \Big),\label{eqn:e308}\\
\kappa^* & = \kappa\Big( \frac{2+\xi}{1-\xi}, \ddot{l}_p(\bgg^*) \Big),
\quad{\rm and}\quad
\tau^* = \frac{9s_0\lambda C_W}{4(1-\xi)(\kappa^*)^2}\quad
{\rm for\ } \xi \in (0,1).
\label{eqn:e309}
\end{align}
Note that $C_W$ is bounded from above.
We closely look at $RE^*$ and $\kappa^*$ in Proposition \ref{prop:prop33}.
Let $\eta^*$ be the smaller solution of
\[
\eta \exp( - \eta ) = \tau^*
\]
as in \cite{HSYYZ2013}. Note that $\tau^*$ should tend to 0 as in Remark 1.

Recall that we are considering $Q_h(\bm{\gamma}; \lambda )$ now since
we can deal with $Q_1(\bm{\gamma}; \lambda )$ in (\ref{eqn:e219}) and
$Q(\bm{\gamma}; \lambda )$ in (\ref{eqn:e231}) in almost the same way and
drive the same results  with just conformable changes.

%%%%%%%%%%%%%%%%%%%%%%%%%%%%%%%%%%%%%%%%%%%%%%%%%% Thm1

\begin{thm}\label{thm:thm31} Assume that Assumptions G and M hold.
Then if $D_l \le \xi \lambda$ for some $\xi \in (0,1)$, we have
\[
P_1( \widehat{\bgg}- \bgg^*) \le \eta^*/C_W.
\]
Then we also have
\begin{align*}
\max_{1\le j \le p}|\widehat g_{cj} - g_{cj}|
& \le C_c \Big( \frac{\eta^*}{L^{1/2}} + L^{-2}\Big),\quad
\max_{1\le j \le p}\|\widehat g_{nj} - g_{nj}\|
\le C_{n1} \Big( \frac{\eta^*}{L^{1/2}} + L^{-2}\Big),\\
\max_{1\le j \le p}\|\widehat g_{nj} - g_{nj}\|_\infty
& \le C_{n2} \Big( \frac{\eta^*}{L^{1/2}} + L^{-2}\Big),
\end{align*}
where $C_c$, $C_{n1}$, and $C_{n2}$ depend on $C_W$, $C_g$, and the properties
of the B-spline basis on $[0,1]$ and they are bounded.
\end{thm}

Some remarks are in order.
%%%%%%%%%%%%%%%%%%%%%%%%%%%% Remark 1

\begin{rmk}
When $p=O(n^{c_p})$ for some $c_p$, we have $D_l=O_p( (n^{-1}\log n)^{1/2})$
and should take $\lambda = C (n^{-1}\log n)^{1/2}$ for some sufficiently
large $C$. As in shown in Proposition \ref{prop:prop33}, we usually have
$(\kappa^*)^2 \sim L^{-1}$ with probability tending to 1 in suitable setups.
Then $\tau^* \sim L (n^{-1}\log n)^{1/2}$ and $\eta^*/\tau^* \to 1$.
This leads to the convergence rate of $O(n^{-2/5}(\log n)^{1/2})$
for $\widehat g_{cj}$ and $\widehat g_{nj}$ and improves
that of \cite{BS2015}, which is $O(n^{-7/20}(\log n)^{1/2})$
for their additive model in a similar setup. Our rate is optimal except for
$(\log n)^{1/2}$. Our results can deal with ultra high-dimensional cases
if $p = \exp ( c_p n )$ and $c_p$ is sufficiently small.
See Propositions \ref{prop:prop32} and \ref{prop:prop33}.
\end{rmk}

%%%%%%%%%%%%%%%%%%%%%%%%%%%% Remark 2
\begin{rmk} Suppose that
\[
\max_{j\in {\cal S}_c}| g_{cj}|/(n^{-2/5}(\log n)^{1/2})\to \infty
\quad{\rm and} \quad
\max_{j\in {\cal S}_n}\| g_{nj}\|/(n^{-2/5}(\log n)^{1/2})\to \infty.
\]
Then if we take $t_\lambda$ satisfying $t_\lambda/\lambda \to \infty$
sufficiently
slowly for $\lambda$ in Remark 1, $\widehat{\cal S}_c$ and $\widehat{\cal S}_n$
in (\ref{eqn:e221}) are consistent estimators of
${\cal S}_c$ and ${\cal S}_n$, respectively.
\end{rmk}

Next we evaluate $D_l$ in Proposition \ref{prop:prop32},
which is called the deviation condition. From Assumption M
and application of Bernstein's inequality
(for example, see \cite{VW1996}), we have
with probability larger than $1-P_Y$,
\begin{equation}
\frac{1}{n}\sum_{i=1}^nY_i(1)=\overline{Y}(1) > C_Y,
\label{eqn:e313}
\end{equation}
where
\[
P_Y= \exp \Big\{ - \frac{C_Y^2n}{2(1+2C_Y/3)} \Big\}.
\]

Since
\begin{equation}
\dot{l}_p( \bgg^* )
= -\frac{1}{n}\sum_{i=1}^n
\int_0^1 \{ \bwi- \widetilde{\bm{W}}_n(t, \bgg^*) \}
dN_i (t),
\label{eqn:e315}
\end{equation}
we evaluate $\dot{l}_{op}$ in (\ref{eqn:e317}) and
$\dot{l}_{op} - \dot{l}_p(\bgg^* )$ in (\ref{eqn:e319}).
\begin{align}
\dot{l}_{op} & =
-\frac{1}{n}\sum_{i=1}^n
\int_0^1\Big\{ \bwi - \frac{S^{(1)}_0(t)}{S^{(0)}_0(t)} \Big\}dN_i(t)
\label{eqn:e317}\\
&= -\frac{1}{n}\sum_{i=1}^n
\int_0^1\Big\{ \bwi - \frac{S^{(1)}_0(t)}{S^{(0)}_0(t)} \Big\}dM_i(t),
\nonumber
\end{align}
where
\[
S^{(k)}_0(t) =\frac{1}{n}\sum_{i=1}^n Y_i(t)
\bm{W}_i^{\otimes k}(t) \exp \{ \bm{g}^T(t) \bxi \}, \quad
{k=0,1,2}.
\]
\begin{equation}
\dot{l}_{op} - \dot{l}_p(\bgg^* )
= \int_0^1\Big\{  \widetilde{\bm{W}}_n(t, \bgg^*)
- \frac{S^{(1)}_0(t)}{S^{(0)}_0(t)} \Big\}d\overline{N}(t).
\label{eqn:e319}
\end{equation}

By combining evaluations of (\ref{eqn:e317}) and (\ref{eqn:e319}),
we obtain Proposition \ref{prop:prop32}. The proof is postponed to
section \ref{sec:proofs}. Recall that $\widetilde{\bm{W}}_n(t, \bgg^*)$
is defined in (\ref{eqn:e301}).

%%%%%%%%%%%%%%%%%%%%%%%%%%%%%%%%%%%%%%%%%%%%%%%%%% Prop.2
\begin{prop}\label{prop:prop32} Assume that Assumptions G and M hold.
Then we have
\[
P_\infty ( \dot{l}_p( \bgg^* ) ) \le
\frac{a_1}{L^{5/2}} + \frac{x(\log n)^{1/2}}{\sqrt{n}}
\]
with probability larger than
\[
1-P_Y - La_2 \exp\{-a_3nL^{-1} \} -2pL \exp
\Big\{ - \frac{a_4x^2\log n }{ 1+x(n^{-1}L\log n)^{1/2} } \Big\},
\]
where $a_j$, $j=1,\ldots,4$, are positive constants depending only on
the assumptions and they are independent of $n$.
\end{prop}

Finally we deal with $\kappa^*$ and $RE^*$.
In Proposition \ref{prop:prop33}, we give their lower bounds.
They are called the compatibility condition and the
restricted eigenvalue condition, respectively.

%%%%%%%%%%%%%%%%%%%%%%%%%%%%%%%%%%%%%%%%%%%%%%%%%% Prop.3
\begin{prop}\label{prop:prop33} Assume that Assumptions G and M hold.
Then with probability larger than $1-P_Y-P_A-P_B-P_C$,
we have
\begin{align*}
\kappa^2( \zeta, \ddot{l}_p( \bgg ) ) & \ge
\exp ( -C_XC_g)(1+O(L^{-2}))\kappa^2( \zeta, \RE \{ \overline{\Sigma } \} )\\
& \qquad - s_0 (1+\zeta )^2 L\Big\{ \frac{c_1}{L^3} + \frac{x(\log n)^{1/2}}{\sqrt{nL}}
\Big\}\\
\intertext{and}
RE^2( \zeta, \ddot{l}_p( \bgg ) )& \ge
\exp ( -C_XC_g)(1+O(L^{-2}))RE^2( \zeta, \RE \{ \overline{\Sigma } \} )\\
& \qquad - s_0 (1+\zeta )^2 L\Big\{ \frac{c_2}{L^3} + \frac{x(\log n)^{1/2}}{\sqrt{nL}}
\Big\}
\end{align*}
where
\begin{align*}
\overline{\Sigma} &= \int_0^1 \overline{G}_Y(t) \lambda_0(t)dt, \quad
\overline{G}_Y(t) = \frac{1}{n}\sum_{i=1}^nY_i(t) \{ \bwi - 
\bm{\mu}_Y(t) \}^{\otimes 2}, \\
\bm{\mu}_Y(t) & = \frac{\RE \{ Y_1(t)\bm{W}_1(t) \}}{\RE \{ Y_1(t) \}},
\quad P_A =2(pL)^2 \exp \Big\{ -\frac{c_3x^2 \log n }{1+ 
x(\log n)^{1/2}(n^{-1}L)^{1/2}} \Big\}, \\
P_B & = 5(pL)^2 \exp \Big\{ -\frac{c_4x (n \log n)^{1/2} }{1+ 
x^{1/2}(n^{-1}\log n)^{1/4}} \Big\},\\
P_C & =2(pL)^2 \exp \Big\{ -\frac{c_5x^2 \log n }{1+ 
x(\log n)^{1/2}n^{-1}} \Big\} .
\end{align*}
Note that $c_j$, $j=1,\ldots, 5$, are positive constants
depending only on the assumptions and they are independent of $n$.
\end{prop}

In the literature, it is often assumed that there
is a positive constant $C_1$ such that
$\lambda_{\min}( \RE \{ \overline{\Sigma} \} ) \ge C_1/L$
due to (\ref{eqn:e901}) and (\ref{eqn:e903}) in \ref{sec:basis}.
Then for some positive constants $C_2$ and $C_3$, we have
\[
\kappa^2( \zeta, \RE \{ \overline{\Sigma } \} ) \ge \frac{C_2}{L}+o_p(L^{-1})
\quad {\rm and} \quad
RE^2 ( \zeta, \RE \{ \overline{\Sigma } \} ) \ge \frac{C_3}{L}+o_p(L^{-1})
\]
if $s_0$ is bounded and $p=O(n^{c_p})$.

%\smallskip
%\medskip
%\begin{equation}\label{eqn:e3}\end{equation}

\section{Other models}
\label{sec:other}
\subsection{Varying coefficient models with index variables}
When we observe $(Z_i(t), \bxi )$ and $Z_i(t)$ is an influential
variable treated as the index variable, the following model
for the compensator is among candidates of our models for
statistical analysis.
\begin{equation}
d\Lambda_i(t) = Y_i(t) \exp \{ g_0(Z_i(t)) + \bxi^T\bm{g}(Z_i(t)) \}
\lambda_0 (t) dt,
\label{eqn:e401}
\end{equation}
where
$Z_i(t) \in [0,1]$, $\int_0^1 g_0(z) dz =0$, and
$g_j(z) =g_{cj} + g_{nj}(z)$, $j=1,\ldots, p$, as in
section \ref{sec:group}. Then we can proceed in almost the
same way with
\begin{align*}
\bwi & = ( \bm{B}^T(Z_i(t)), \bm{X}_i(t)^T\otimes
\overline{\bm{B}}^T( Z_i(t) ))^T,\\
\bgg & = ( \bgg_{-10}^T, \gamma_{11},  \bgg_{-11}^T,
\ldots, \gamma_{1p},  \bgg_{-1p}^T)^T,\\
P_1(\bm{\gamma}) & = \sum_{j=0}^p |\gamma_{1j}| + \sum_{j=1}^p | \bgmj |,\\
P_h(\bm{\gamma}) & = \sum_{j=1}^p ( |\gamma_{1j}|^q + | \bgmj |^q)^{1/q}
+  \sum_{j=0}^p | \bgmj | ,\\
P_\infty (\bm{\gamma}) & = \{\max_{1\le j \le p}|\gamma_{1j}| \vee |\bgmj| 
\}\vee | \bgg_{-10} |, \\
Q_1(\bgg ;\lambda ) & = l_p (\bgg ) + \lambda P_1( \bgg ), \quad
{\rm and}\quad Q_h(\bgg ;\lambda ) = l_p (\bgg ) + \lambda P_h( \bgg ).
\end{align*}

We can carry out simultaneous variable selection and structure identification
of this model as for time-varying coefficient models and we are able
to prove the same results in almost the same way. Almost no
change is necessary to the proofs of Proposition
\ref{prop:prop31} and Theorem \ref{thm:thm31}. When we consider
Propositions \ref{prop:prop32} and \ref{prop:prop33}, 
we should be a little careful in evaluating predictable variation processes
and so on.
Then we have to deal with terms like
\[
n^{-1}\sum_{i=1}^n |b_{0j}(Z_i(t))|,\quad
n^{-1}\sum_{i=1}^n |b_{j}(Z_i(t))|,\quad{\rm and}\quad
n^{-1}\sum_{i=1}^n |b_{j}(Z_i(t))b_{k}(Z_i(t))|
\]
as compared to
\[ |b_{0j}(t)|, \quad  |b_{j}(t)|,  
\quad{\rm and}\quad |b_{j}(t)b_{k}(t)|
\]
for time-varying coefficient models. Note that
we can use exponential inequalities for generalized U-statistics as
given in Gine et al.\cite{GLZ2000} instead of Lemma 4.2 in \cite{HSYYZ2013}
in the proof of Proposition \ref{prop:prop33}. We give more details in
\ref{sec:additionalproofs}.

\subsection{Additive models}
When we have no specific index variable,
the following additive model may be suitable. 
\begin{equation}
d\Lambda_i(t) = Y_i(t) \exp \Big\{ \sum_{j=1}^p g_j(X_{ij}(t)) \Big\}
\lambda_0 (t) dt,
\label{eqn:e403}
\end{equation}
where $\int_0^1 g_j(x) dx =0$ and $X_{ij}(t) \in [0,1]$.
These $g_j(x)$ can be orthogonally decomposed into the
linear part and the nolinear part as well.

We should take $b_2(X_{ij}(t)) = (12L^{-1})^{1/2}(X_{ij}(t)-1/2)$ and use
$b_2(X_{ij}(t))$ and $( b_3(X_{ij}(t)), \ldots, b_L(X_{ij}(t)))^T$ for the linear part and the nonlinear
part, respectively. We have no $b_1(X_{ij}(t)) $ and divide $\bgg_{-1j}$ into
$\gamma_{2j}$ and $\bm{\gamma}_{-2j}=(\gamma_{3j},\ldots, \gamma_{Lj})^T$.
Then we can apply the same group Lasso procedure for
variable selection and structure identification with
\begin{align*}
\bwi & = ( \bm{B}^T(X_{i1}(t)),\ldots, \bm{B}^T(X_{ip}(t)))^T,\quad
\bgg_{-1}  = ( \bgg_{-11}^T, \ldots, \bgg_{-1p}^T)^T,\\
P_1(\bgg_{-1}) & = \sum_{j=1}^p | \gamma_{2j}| +\sum_{j=1}^p | \bm{\gamma}_{-2j} |,\\
P_h(\bgg_{-1}) & = \sum_{j=1}^p ( |\gamma_{2j}|^q + | \bm{\gamma}_{-2j} |^q)^{1/q}
+  \sum_{j=1}^p | \bm{\gamma}_{-2j}  |, \\
P_\infty (\bgg_{-1} ) & = \max_{1\le j \le p}  |\gamma_{2j}|
\vee | \bgg_{-2j} |,\\
Q_1(\bgg_{-1};\lambda ) & = l_p (\bgg_{-1} ) + \lambda P_1( \bgg_{-1} ),
\quad {\rm and}\quad 
Q_h(\bgg_{-1};\lambda ) = l_p (\bgg_{-1} ) + \lambda P_h( \bgg_{-1} ).
\end{align*}

We have the same theoretical results with just conformable changes.
We should be careful in the proofs of Propositions \ref{prop:prop32}
and \ref{prop:prop33} as for varying coefficient models with index variables,
too. We have to deal with terms like
\[
n^{-1}\sum_{i=1}^n |b_{0j}(X_{il}(t))|,\quad
n^{-1}\sum_{i=1}^n |b_{j}(X_{il}(t))|,\quad{\rm and}\quad
n^{-1}\sum_{i=1}^n |b_{j}(X_{il}(t))b_{k}(X_{il}(t))|
\]
as compared to
\[ |b_{0j}(t)|, \quad  |b_{j}(t)|,  
\quad{\rm and}\quad |b_{j}(t)b_{k}(t)|
\]
for time-varying coefficient models. We can use exponential inequalities
for generalized U-statistics as given in Gine et al.\cite{GLZ2000} instead of Lemma 4.2
in \cite{HSYYZ2013} in the proof of Proposition \ref{prop:prop33}.

%\smallskip
%\medskip
%\begin{equation}\label{eqn:e3}\end{equation}

\section{Numerical studies}
\label{sec:numerical}
%!TEX root = Coxlasso.tex
We carried out a small simulation study for the two models in
section \ref{sec:other} with the $P_1$ penalty because
time-varying coefficient models and the $P_h$ penalty are
numerically intractable at present and our computational ability
is limited. We used the grpsurv function of the package ``grpreg"
(Breheny\cite{Brehney2016}) for R in our numerical study and
all the covariates are time-independent. An extensive numerical study
is a topic of future research. 

First we describe the data generating process of the covariates :
$\{ X_{ij} \}_{j=1}^q$, $\{ X_{ij} \}_{j=q+1}^p$, and $Z_i$
are mutually independent. Then
%$\{ X_{ij} \}_{j=q+1}^p$ 
$X_{ij}$, $j=q+1, \ldots, p$, and $Z_i$
follow $U(0,1)$ independently. We define $\{ X_{ij} \}_{j=1}^q$
in (\ref{eqn:e501}).
\begin{equation}
X_{ij}=F(Y_{ij}),\quad j=1,\ldots, q,
\label{eqn:e501}
\end{equation}
where $\{Y_{ij} \}$ is a stationary Gaussian AR(1) process
with $\rho =0.3$ and $F(y)$ is the distribution function of
$Y_{ij}$.

%%%%%%%%%%%%%%%%%%% VC model
Next we gives the details for our varying coefficient model
with an index variable $Z$. We took
%1.VC 関数、lambda_0、C_i
\[
\lambda_0(t) = 0.5, \quad g_1(z)=g_2(z)=1, \quad
g_3(z) =4z, \quad g_4(z) =4z^2. 
\]
The other functions are taken to be $0$. Hence we have $s_c=4$ and $s_n =2$.
Note that $X_1$ and $X_2$ are relevant for only the constant component and that
$X_3$ and $X_4$ are relevant for both the constant component and
the non-constant one. All the other covariates are irrelevant. We imposed no penalty on
the coefficient vector for $g_0(z)$ in this simulation study. The censoring variable
$C_i$ follows the exponential distribution with mean$ =1/0.85$ independently of all
the other variables and the censoring rate is about $ 20\%$. 

%%%%%%%%%%%%%%%%%%% Additive model
Then we describe the details for our additive model. We took
%2.Additive 関数、lambda_0、C_i
\begin{align*}
%\[
\lambda_0(t) &= 0.5,\quad g_1(x) = g_2(x)=2^{1/2}(x-1/2),\\
g_3(x) & = 2^{-1/2}\cos (2\pi x) + (x-1/2),\quad 
g_4(x) = \sin (2\pi x). 
%\]
\end{align*}
The other functions are taken to be $0$. Hence we have $s_c=4$ and $s_n =2$
and note that $X_1$ and $X_2$ are relevant for only the linear component and
that $X_3$ and $X_4$ are relevant for both the linear component and
the nonlinear one. All the other covariates are irrelevant. The censoring
variable $C_i$ follows the exponential distribution with mean$ =1/0.80$
independently of all the other variables and the censoring rate is about $ 30 \%$. 

%%%%%%%%%%%%%%%%%% Parameters & 表の説明 & no tuning parameter selection
When we carried out simulations,
we took $p=400$, $q=8$, and $L=6$. We used the quadratic spline basis
and the repetition number is $500$. The results are given in Tables 1 and 2.
In addition to the group Lasso, we applied a threshold method in (\ref{eqn:e221})
with $t_\lambda =0.1$. In the tables, $t_\lambda =0$ means the group Lasso
and $t_\lambda =0.1$ means this threshold group Lasso. In the tables, Failure,
Correct, and False respectively stand for 

\smallskip
{\bf Failure}: The rate of relevant covariates that are not chosen wrongly,

{\bf Correct}: The rate of correct decisions,

{\bf False}: The rate of irrelevant covariates that are wrongly chosen.

\smallskip
As for the tuning parameter $\lambda$, we tried several values and found
variable selection and structure identification are sensitive to this
$\lambda$. We presented one of the good results for each model here.
In Table 2, we sometimes missed the linear components of $X_3$ and $X_4$.
If we incorporate the assumption in (\ref{eqn:e205}), we will not 
miss these linear components. Since our procedure
can be seen as a screening procedure, screening consistency
or not to miss any relevant covariates is inevitable.
When $p$ is very large compared to $n$, it may be better
to consider only variable selection based
on (\ref{eqn:e231}) first and then apply
our procedure based on some weighted $P_1(\bm{\gamma})$
as in the adaptive group Lasso.

%%%%%%%%%%%%%%%%%%% Tuning parameter \lambda
As for tuning parameter selection rules,
we don't have any results on them at present although the results of Tables 1
and 2 seem to be very promising. Some rules based on BIC, the number of selected
variables, analysis of solution paths, a threshold value method, or combinations
of them may be possible for screening consistency, not for selection consistency.
These rules are a topic of future research since our orthonormal basis
method of simultaneous variable selection and structure identification for (ultra)
high-dimensional Cox models has just been proposed.

\medskip
\begin{table}[ht]
\begin{flushleft}
\scalebox{0.8}{
%\begin{tabular}{lllllllll}\hline 
 \begin{tabular}{lcccccccc}\hline 
$\lambda=0.08$&\multicolumn{2}{c}{$X_1$ and $X_2$}
&\multicolumn{2}{c}{$X_3$ and $X_4$}& \multicolumn{2}{c}{
$X_{5}$ to $X_q(q=8)$}&\multicolumn{2}{c}{$X_{q+1}$ to $X_p(p=400)$}\\ 
\hline
 $t_\lambda = 0$& Const. & Non-const. & Const. & Non-const. & Const. & Non-const. & Const. & Non-const.\\ 
  \hline
Failure & 0.000 & --- & 0.000 & 0.016 & --- & --- & --- & --- \\ 
Correct & 1.000 & 0.993 & 1.000 & 0.984 & 0.948 & 0.988 & 0.954 & 0.996 \\ 
False & --- & 0.007 & --- & --- & 0.052 & 0.012 & 0.046 & 0.004 \\ 
   \hline
   $t_\lambda = 0.1$ &Const. & Non-const. & Const. & Non-const. & Const. & Non-const. & Const. & Non-const. \\ 
  \hline
Failure & 0.001 & --- & 0.000 & 0.029 & --- & --- & --- & --- \\ 
Correct & 0.999 & 0.996 & 1.000 & 0.971 & 0.968 & 0.997 & 0.974 & 0.998 \\ 
False & --- & 0.004 & --- & --- & 0.032 & 0.003 & 0.026 & 0.002 \\ 
   \hline   \end{tabular}
}
\caption{Varying coefficient model with an index variable} 
\label{SimuVC}\end{flushleft}
\end{table}

\begin{table}[!h]
\begin{flushleft}
\scalebox{0.8}{
%\begin{tabular}{lllllllll}\hline
 \begin{tabular}{lcccccccc}\hline 
$\lambda=0.1$&\multicolumn{2}{c}{$X_1$ and $X_2$}&\multicolumn{2}{c}{$X_3$ and $X_4$}&
\multicolumn{2}{c}{$X_{5}$ to $X_q(q=8)$
}&\multicolumn{2}{c}{$X_{q+1}$ to $X_p(p=400)$}\\    \hline
$t_\lambda=0$ & Linear & Nonlinear & Linear & Nonlinear & Linear & Nonlinear & Linear & Nonlinear \\ 
  \hline
Failure & 0.000 & --- & 0.065 & 0.000 & --- & --- & --- & --- \\ 
Correct & 1.000 & 0.900 & 0.935 & 1.000 & 0.994 & 0.932 & 0.997 & 0.926 \\ 
False & --- & 0.100 & --- & --- & 0.006 & 0.068 & 0.003 & 0.074 \\ 
    \hline
 $t_\lambda =0.1$& Linear & Nonlinear & Linear & Nonlinear & Linear & Nonlinear & Linear & Nonlinear \\ 
  \hline
Failure & 0.000 & --- & 0.181 & 0.000 & --- & --- & --- & --- \\ 
Correct & 1.000 & 0.983 & 0.819 & 1.000 & 1.000 & 0.992 & 1.000 & 0.987 \\ 
False & --- & 0.017 & --- & --- & 0.000 & 0.008 & 0.000 & 0.013 \\ 
   \hline
\end{tabular}
}
\caption{Additive model}
\label{SimuAdditive}

\end{flushleft}
\end{table}

\section{Proofs}
\label{sec:proofs}
We prove Propositions \ref{prop:prop31}-\ref{prop:prop33}
and Theorem \ref{thm:thm31}.

For a vector $\bm{a}$ and a matrix $A$,
$(\bm{a})_i$ and $(A)_{ij}$ mean the $i$th element
of $\bm{a}$ and the $(i,j)$ element of $A$, respectively.
We present the proofs of technical lemmas in \ref{sec:technical}.

%%%%%%%%%%%%%%%%%%%%%%% Proof of Prop 31
\begin{pop1} Note that
\begin{eqnarray}
\lefteqn{ (\widehat{\bgg} - \bgg^*)^T (\dot{l}_p (\widehat{\bgg})
- \dot{l}_p(\bgg^*)} \label{eqn:e601}\\
& = & \Big\{  \sum_{j \in \overline{{\cal S}}_c }\wt_{1j}
\frac{\partial l_p}{\partial \gamma_{1j}} ( \widehat{\bgg} )
+ \sum_{j \in \overline{{\cal S}}_c } \wbt_{-1j}^T
\frac{\partial l_p}{\partial \bgmj }( \widehat{\bgg} ) \Big\}
\nonumber\\
&  & + \Big\{  \sum_{j \in \overline{{\cal S}}_n \cap {\cal S}_c }\wt_{1j}
\frac{\partial l_p}{\partial \gamma_{1j}}(\widehat{\bgg})
+  \sum_{j \in \overline{{\cal S}}_n\cap {\cal S}_c }\wbt_{-1j}^T
\frac{\partial l_p}{\partial \bgmj }(\widehat{\bgg}) \Big\}
\nonumber\\
&  & + \Big\{  \sum_{j \in {\cal S}_n }\wt_{1j}
\frac{\partial l_p}{\partial \gamma_{1j}}(\widehat{\bgg})
+  \sum_{j \in {\cal S}_n }\wbt_{-1j}^T
\frac{\partial l_p}{\partial \bgmj }(\widehat{\bgg}) \Big\}
\nonumber\\
& & + \{ - \wbt^T(\dot{l}_p( \bgg^* ) \}
= E_1 + E_2 + E_3 + E_4 \ge 0. \quad {\rm (say)}
\nonumber
\end{eqnarray}
The last inequality follows from the convexity of $l_p( \bgg )$
and we should recall that $\wbt = \widehat{\bgg} -\bgg^*$.

We evaluate $E_j$, $j=1,2,3,4$.

%%%%%%%%%%%%%%%%%%%%%% E_1
\noindent
${\bf E_1:}$ Notice that $\widehat{\bgg}_j = \wbt_j$. Then we should
evaluate 
\[
E_{1j}= \wt_{1j}
\frac{\partial l_p}{\partial \gamma_{1j}} ( \widehat{\bgg} )
+  \wbt_{-1j}^T
\frac{\partial l_p}{\partial \bgmj }( \widehat{\bgg} ).
\]
Recalling (\ref{eqn:e220}), we use the results in \ref{sec:subgradients}.

When $\widehat{\gamma}_{1j}\ne 0$ and $\widehat{\bgg}_{-1j}\ne 0$, we have
\begin{equation}
E_{1j} = -\lambda ( |\wt_{1j} |^q +  |\wbt_{-1j} |^q )^{1/q} -
\lambda | \wbt_{-1j}| .
\label{eqn:e603}
\end{equation}
When $\widehat{\gamma}_{1j}\ne 0$ and $\widehat{\bgg}_{-1j}= 0$, we have
\begin{equation}
E_{1j}= - \lambda |\wt_{1j} | .
\label{eqn:e605}
\end{equation}
When $\widehat{\gamma}_{1j}= 0$ and $\widehat{\bgg}_{-1j}\ne 0$, we have
\begin{equation}
E_{1j}= - 2 \lambda |\wbt_{-1j} | .
\label{eqn:e607}
\end{equation}
From (\ref{eqn:e603})-(\ref{eqn:e607}), we obtain
\begin{equation}
E_1 \le -\lambda \sum_{j \in \overline{{\cal S}}_c }
( |\wt_{1j} | +  |\wbt_{-1j} | ) .
\label{eqn:e609}
\end{equation}

%%%%%%%%%%%%%%%%%%%%%% E_2
\noindent
${\bf E_2:}$ Notice that $\widehat{\bgg}_{-1j} = \wbt_{-1j}$
and $| \frac{\partial l_p}{\partial \gamma_{1j}} ( \widehat{\bgg} ) | 
\le \lambda $. Then we should evaluate 
\[
E_{2j}= \wt_{1j}
\frac{\partial l_p}{\partial \gamma_{1j}} ( \widehat{\bgg} )
+  \widehat{\bgg}_{-1j}^T
\frac{\partial l_p}{\partial \bgmj }( \widehat{\bgg} ).
\]

When $ \widehat{\gamma}_{1j} \ne 0$ and $ \widehat{\bgg}_{-1j} \ne 0$,
we have
\begin{align}
E_{2j} & \le \lambda |\wt_{1j} | -\lambda
( | \widehat{\gamma}_{ij}|^q  + |\wbt_{-1j} |^q )^{\frac{1}{q}-1} | \wbt_{-1j}|^q
- \lambda | \wbt_{-1j}| \label{eqn:e611}\\
& \le \lambda ( | \wt_{1j} | - | \wbt_{-1j} |) .\nonumber
\end{align}
When $ \widehat{\gamma}_{1j} \ne 0$ and
$ \widehat{\bgg}_{-1j}= 0$, we have
\begin{equation}
E_{2j}\le  \lambda |\wt_{1j} | .
\label{eqn:e613}
\end{equation}
When $ \widehat{\gamma}_{1j}= 0$ and $ \widehat{\bgg}_{-1j}\ne 0$
and when $ \widehat{\gamma}_{1j}= 0$ and $ \widehat{\bgg}_{-1j}= 0$,
we have
\begin{equation}
E_{2j}\le \lambda |\wt_{1j} | - 2 \lambda |\wbt_{-1j} | .
\label{eqn:e615}
\end{equation}
From (\ref{eqn:e611})-(\ref{eqn:e615}), we obtain
\begin{equation}
E_2 \le \lambda \sum_{ j \in \overline{{\cal S}}_n \cap {\cal S}_c
 } ( |\wt_{1j} | -  |\wbt_{-1j} | )
\le  \lambda \sum_{ j \in \overline{{\cal S}}_n \cap {\cal S}_c
 } ( 2|\wt_{1j} | -  |\wbt_{-1j} | ) .
\label{eqn:e617}
\end{equation}

%%%%%%%%%%%%%%%%%%%%%% E_3
\noindent
${\bf E_3:}$ Notice that 
$| \frac{\partial l_p}{\partial \gamma_{1j}} ( \widehat{\bgg} ) | 
\le \lambda $ and
$ |\frac{\partial l_p}{\partial \bgmj }( \widehat{\bgg} )| \le 2 \lambda $.
Then we have
\begin{equation}
E_3 \le 2\lambda \sum_{ j \in {\cal S}_n } ( |\wt_{1j} | + |\wbt_{-1j} | ). 
\label{eqn:e619}
\end{equation}

%%%%%%%%%%%%%%%%%%%%%% E_4
\noindent
${\bf E_4:}$ We have
\begin{equation}
E_4 \le P_1 (\wbt) D_l = ( P_1(\wbt_{\cal S}) + P_1 (\wbt_{\overline{{\cal S}}})
)D_l. \label{eqn:e621}
\end{equation}

(\ref{eqn:e609}), (\ref{eqn:e617}), (\ref{eqn:e619}), and (\ref{eqn:e621})
yield that
\[
E_1 + E_2 + E_3 + E_4 \le (2\lambda + D_l) P_1(\wbt_{\cal S})
-(\lambda - D_l) P_1 (\wbt_{\overline{{\cal S}}}).
\]
The first and second inequalities follow from (\ref{eqn:e601})
and the above inequality. The third inequality
follows from the following expression of the second one.
\[
P_1 (\wbt_{\overline{{\cal S}}}) \le \frac{2\lambda + D_l}{\lambda - D_l}
P_1 (\wbt_{\cal S})
\]
Hence the proof of the proposition is complete.
\end{pop1}

We establish the oracle inequality.
%%%%%%%%%%%%%%%%%%%%%%%%%%%%%%%%% Thm 1

\begin{pot} First we define $D(\wt )$ by
\[
D( \bt ) = \max_{i,j}\max_{0\le t \le 1}| \bt^T \bwi - \bt^T \bm{W}_j(t)|. 
\]

We need two lemmas.
\begin{lem}\label{lem:lem71}
\[
D( \bt ) \le C_W P_1 (\bt )
\]
\end{lem}
\begin{lem}\label{lem:lem72}
\[
e^{ -D (\bt ) } \bt^T\ddot{l}_p(\bgg^* )\bt
\le (\bgg^* + \bt  - \bgg^*)^T (\dot{l}_p (\bgg^* + \bt)
- \dot{l}_p(\bgg^*)
\le e^{D (\bt )} \bt^T\ddot{l}_p(\bgg^* )\bt
\]
\end{lem}

Now we begin to prove the oracle inequality. If $\wbt =0$,
the desired inequality holds. Hence we assume $\wbt \ne0$
and set
\[
\wbb = \frac{\wbt}{P_1(\wbt)}.
\]

We have from Proposition \ref{prop:prop31} and the definition of
$P_1(\bgg )$ that
\begin{equation}
\wbb \in \Theta \Big( \frac{2+\xi}{1-\xi} \Big)
\quad {\rm and} \quad P_1( \wbb) =
P_1(\wbb_{\cal S}) +  P_1(\wbb_{\overline{\cal S}})=1.
\label{eqn:e623}
\end{equation}

When $D_l \le \xi \lambda$, the first inequality of Proposition
\ref{prop:prop31} implies that the following inequalities hold
at $x=0$ and $x=P_1(\wbt )$.
\begin{eqnarray}
\lefteqn{\wbb^T \{ \dot{l}(\bgg^* + x \wbb) - \dot{l}(\bgg^*)\}
}\label{eqn:e625}\\
& \le & (2+\xi ) \lambda P_1 (\wbb_{\cal S}) - (1-\xi ) \lambda
P_1(\wbb_{\overline{\cal S}})
\nonumber\\
& = & 3\lambda P_1 (\wbb_{\cal S}) - \lambda(1-\xi) \,
\le\, \frac{9\lambda}{4(1-\xi)} \{ P_1( \wbb_{\cal S}) \}^2.
\label{eqn:e627}
\end{eqnarray}
We also used (\ref{eqn:e623}) here.

Note that (\ref{eqn:e625}) is monotone increasing and continuous
in $x$ due to the convexity of $l_p( \bgg )$ and we have (\ref{eqn:e627})
on $[0, P_1(\wbt )]$. Let $x_{\bm{b}}$ be the maximum of $x$ satisfying
\begin{equation}
\wbb^T \{ \dot{l}(\bgg^* + x \wbb) - \dot{l}(\bgg^*)\}
\le \frac{9\lambda}{4(1-\xi)} \{ P_1( \wbb_{\cal S}) \}^2
\label{eqn:e629}
\end{equation}
for any $s\in [0, x]$.

If we find an upper bound of $x_{\bm{b}}$, say $x_0$,
we have $P_1(\wbt ) \le x_0$. Therefore we will find
an upper bound of $x_{\bm{b}}$ as in \cite{HSYYZ2013}.

From Lemmas \ref{lem:lem71} and \ref{lem:lem72}, we have
\begin{align}
x \wbb^T \{ \dot{l}(\bgg^* + x \wbb) - \dot{l}(\bgg^*)\}
& \ge x^2 \exp \{ -D(x\wbb ) \} \wbb^T \ddot{l}_p (\bgg^* )
\wbb
\label{eqn:e639}\\
& \ge x^2 \exp \{ -C_Wx \} \wbb^T \ddot{l}_p (\bgg^* )
\wbb .\nonumber
\end{align}

The definition of $\kappa^*$ and (\ref{eqn:e639}) imply that
\begin{equation}
\wbb^T \{ \dot{l}(\bgg^* + x \wbb) - \dot{l}(\bgg^*)\}
\ge x \exp \{ -C_Wx \} \frac{(\kappa^*)^2}{s_0}
\{ P_1( \wbb_{\cal S}) \}^2 .
\label{eqn:e641}
\end{equation}

It follows from (\ref{eqn:e627}) and (\ref{eqn:e641}) that
\[
\frac{9\lambda s_0 C_W}{4(1-\xi)(\kappa^*)^2}= \tau^*
\ge C_W x\exp \{ -C_Wx \}.
\]
Consequently we have from the definition of $\eta^*$ and
the above inequality that
\[
C_Wx_{\bm{b}} \le \eta^* \quad {\rm and} \quad
\frac{\tau^*}{\eta^*} \to 1 \ {\rm if}\ \tau^* \to 0.
\]
We have found that $ \eta^*/C_W $ is an upper bound of $x_{\bm{b}}$
and that $P_1(\wbt ) \le \eta^*/C_W$.

As for the the rest of the theorem, the result on $\widehat{g}_{cj}$
is straightforward from (\ref{eqn:e218}). The upper bounds on
$\widehat{g}_{nj}(t)$ follow from (\ref{eqn:e901}),
%(\ref{eqn:e905}),
(\ref{eqn:e904}), and the following inequalities.
\begin{align*} |(\widehat{\bgg}_{-1j}
- \bgg_{-1j}^* )^T \bm{B}(t) |
& \le \{ \lambda_{\max}(A_{-1}A_{-1}^T)\}^{1/2}|\widehat{\bgg}_{-1j}
- \bgg_{-1j}^*| | \bm{B}_0(t)  |\quad
{\rm and}\\ | \bm{B}_0(t) | & \le 1
\end{align*}
Recall that the properties of our basis are collected
in \ref{sec:basis}.

Hence the proof of the theorem is complete.
\end{pot}

Now we prove Proposition \ref{prop:prop32}.

%%%%%%%%%%%%%%%%%%%%%%% Proof of Prop 32
\begin{pop2} We implicitly carry out our evaluation on
$\{ \overline{Y}(1) > C_Y \}$. $C_1,C_2,\ldots$ are
generic positive constants and they depend only on
the assumptions.

First we deal with (\ref{eqn:e319}), which is represented as
\begin{equation}
\int_0^1
\Big[\frac{ S_0^{(0)}(t)\{ S^{(1)}(t,\bgg^*) - S_0^{(1)}(t) \}
}{ S^{(0)}(t,\bgg^*)S_0^{(0)}(t) }+
\frac{ S_0^{(1)}(t) \{ S_0^{(0)}(t)- 
S^{(0)}(t,\bgg^*) \} }{ S^{(0)}(t,\bgg^*)S_0^{(0)}(t) }
\Big]d \overline{N}(t) .
\label{eqn:e643}
\end{equation}
We can rewrite the expression in (\ref{eqn:e643}) as
\begin{align}
(\ref{eqn:e643})& = (I\otimes A_0)
\int_0^1
\Big[\frac{ S_0^{(0)}(t)\{ \overline{S}^{(1)}(t,\bgg^*) - \overline{S}_0^{(1)}(t) \}
}{ S^{(0)}(t,\bgg^*)S_0^{(0)}(t) } \label{eqn:e642}\\
& \quad + 
\frac{ \overline{S}_0^{(1)}(t) \{ S_0^{(0)}(t)- 
S^{(0)}(t,\bgg^*) \} }{ S^{(0)}(t,\bgg^*)S_0^{(0)}(t) }
\Big]d \overline{N}(t)\nonumber\\
&=  (I\otimes A_0) \Delta l_p \quad{\rm (say)},\nonumber
\end{align}
where
\begin{align*}
\overline{S}^{(1)}(t,\bm{\gamma})&
=\frac{1}{n}\sum_{i=1}^n Y_i(t)
(\bm{X}_i(t)\otimes \bm{B}_0(t)) \exp \{ \bm{W}_i^T(t) \bgg \},
\\
\overline{S}_0^{(1)}(t)&= \frac{1}{n}\sum_{i=1}^n Y_i(t)
(\bm{X}_i(t)\otimes \bm{B}_0(t))\exp \{ \bm{X}_i(t)^T\bm{g}(t) \}.
\end{align*}

Due to the definition of $\bgg^*$, we have uniformly in
$t$ and $l(0\le l < p)$,
\[ | S_0^{(0)}(t) - S^{(0)}(t,\bgg^*) | \le C_1 L^{-2},
%\quad {\rm and}\quad 
\ C_2 \le S_0^{(0)}(t)\wedge S^{(0)}(t,\bgg^*),
\ S_0^{(0)}(t)\vee S^{(0)}(t,\bgg^*) \le C_3, 
\]
\[  | ( \overline{S}_0^{(1)}(t) - \overline{S}^{(1)}(t,\bgg^*))_{lL+j} |
\le C_4L^{-2} |b_{0j} (t) |,
\]
\[ |(\overline{S}_0^{(1)}(t))_{lL+j}| \vee | (\overline{S}^{(1)}(t,\bgg^*))_{lL+j}|
\le C_5  |b_{0j}(t) |.
\]
%
%For the $(lL+j)$th element of (\ref{eqn:e643}) with $j=1$,
%it is easy to see that they are uniformly bounded from above
%by $C_6L^{-5/2}$.
%the $(lL+j)$th element of (\ref{eqn:e643}) with $j\ge 2$,

Now we evaluate $\Delta l_p$. Its $(lL+j)$th element is bounded from above
by 
\begin{equation}
C_6L^{-2}\int_0^1 | b_{0j}(t) | d \overline{N} (t).
\label{eqn:e646}
\end{equation}
for some positive constant $C_6$. First notice that
\begin{equation}
\int_0^1  |b_{0j}(t) | d \overline{N} (t) = \int_0^1  |b_{0j}(t) | d \overline{M}(t)
+ O(L^{-1})
\label{eqn:e644}
\end{equation}
uniformly in $j$. Then application of an exponential inequality for martingales
(Lemma 2.1 in \cite{VDG1995}) yields
\begin{equation}
\RP \Big(
\max_{2\le j \le L}\int_0^1  |b_{0j}(t) | d \overline{M}(t) > \frac{x}{L}
\Big)\le
L C_7 \exp \Big\{ -C_8 \frac{nL^{-1}x^2}{1+x} \Big\}.
\label{eqn:e645}
\end{equation}
We used the properties of the support of the B-spline basis in (\ref{eqn:e644})
and (\ref{eqn:e645}). Taking $x=1$ in (\ref{eqn:e645}), we have established
\begin{equation} | \Delta l_p |_\infty
%\max_{1\le j \le L,\ 0\le l <p} | (\Delta l_p)_{lL+j} |
\le \frac{C_9}{L^3}
\label{eqn:e647}
\end{equation}
with probability larger than $1 - L C_7 \exp \Big\{ -2^{-1}C_8 nL^{-1} \Big\}$.

From (\ref{eqn:e642}), (\ref{eqn:e647}), and (\ref{eqn:e906}), we obtain
\begin{equation}
P_\infty ( \dot{l}_{op}- \dot{l}_p( \bgg^* ) ) \le C_{10}L^{-5/2}
\label{eqn:e649}
\end{equation}
with probability larger than $1 - L C_7 \exp \Big\{ -2^{-1}C_8 nL^{-1} \Big\}$.

Finally we deal with (\ref{eqn:e317}) by exploiting the same
exponential inequality for martingales.

For the $(lL+j)$th element with
$j=1$, we have
\begin{equation}
\RP \Big( |(\dot{l}_{op})_{lL+j} | \ge
\frac{x(\log n)^{1/2}}{\sqrt{nL}} \Big)\le 2 \exp \Big\{
-\frac{C_{11}x^2\log n}{x(n^{-1}\log n)^{1/2}+1}
\Big\}.
\label{eqn:e651}
\end{equation}

For the $(lL+j)$th element with
$j\ge 2$, we have
\begin{equation}
\RP \Big( |(\dot{l}_{op})_{lL+j} | \ge
\frac{x(\log n)^{1/2}}{\sqrt{nL}} \Big)\le
2 \exp \Big\{
-\frac{C_{12}x^2\log n}{x(n^{-1}L\log n)^{1/2}+1} \Big\}.
\label{eqn:e653}
\end{equation}
We used the fact that
\begin{equation}
\int_0^1 b_j^2(t) \lambda_0(t) dt \le C_\lambda \bm{a}_{0j}^T
\Omega_0  \bm{a}_{0j} = O(L^{-1})
\label{eqn:e690}
\end{equation}
when we evaluated the predictable variation process.

It follows from (\ref{eqn:e651}) and (\ref{eqn:e653}), 
that
\begin{equation}
P_\infty ( \dot{l}_{op} )  \le  x(\log n)^{1/2}n^{-1/2} 
\label{eqn:e654}
\end{equation}
with probability larger than 
\begin{equation}
1 - 2pL\exp \Big\{
-\frac{C_{13}x^2\log n}{x(n^{-1}L\log n)^{1/2}+1} \Big\}.
\label{eqn:e658}
\end{equation}
Hence the desired result follows from (\ref{eqn:e313}),
(\ref{eqn:e649}), and (\ref{eqn:e654}) and
the proof of the proposition is complete.
\end{pop2}

Finally we give the proof of Proposition \ref{prop:prop33}.

%%%%%%%%%%%%%%%%%%%%%%% Proof of Prop 33
\begin{pop3} $C_1,C_2,\ldots$ are
generic positive constants and they depend only on
the assumptions. We use the following lemma, which is
a version of Lemma 4.1(ii) in \cite{HSYYZ2013}.

\begin{lem}\label{lem:lem73} 
\begin{align*}
\kappa^2 (\zeta,\Sigma_1) & \ge \kappa^2 (\zeta, \Sigma_2) -s_0
(1+\zeta )^2 L \max_{j,k}| ( \Sigma_1 - \Sigma_2 )_{jk} |\\
RE^2 (\zeta,\Sigma_1) & \ge RE^2 (\zeta, \Sigma_2) -s_0
(1+\zeta )^2 L \max_{j,k}| ( \Sigma_1 - \Sigma_2 )_{jk} |
\end{align*}
When $\Sigma_2 - \Sigma_1 $ is n.n.d., we can replace
$\Sigma_2 - \Sigma_1 $ in the above inequalities with
$\Delta$ such that $\Delta - ( \Sigma_2 - \Sigma_1 )$
is n.n.d.
\end{lem}

We implicitly carry out our evaluation on $\{ \overline{Y}(1) > C_Y \}$.
First we outline the proof and then give the details.

Define $\widetilde{\Sigma}_0$ by
\begin{align}
\widetilde{\Sigma}_0 &= \int_0^1 V_n(t, \bgg^* ) S_0^{(0)}(t) \lambda_0 (t) dt
\label{eqn:e655}\\
\intertext{and set}
\Delta_1 & = \ddot{l}_p (\bgg^* ) - \widetilde{\Sigma}_0 =
\int_0^1 V_n(t, \bgg^* ) d \overline{M}(t).
\label{eqn:e656}
\end{align}
We treat $\Delta_1 $ by using the exponential inequalities for
martingales.

Next define $\widetilde{\Sigma}$ by
\[
\widetilde{\Sigma} = \int_0^1 V_n(t, \bgg^* ) S^{(0)}(t,\bgg^*)
 \lambda_0 (t) dt
\]
and set $ \Delta_2 = \widetilde{\Sigma}_0 - \widetilde{\Sigma}$.
Since
\[ |\bm{W}_i^T(t)\bgg^* - \bm{X}_i^T(t) \bgt | \le C_X C_{approx}L^{-2}
\]
and we can use the results on predictable variation process in
evaluating $\Delta_1$, we can easily prove
\begin{equation}
\max_{j,k}| (\Delta_2)_{jk}| \le C_1 L^{-3}.
\label{eqn:e657}
\end{equation}
We omit the details for (\ref{eqn:e657}) in this paper.

Define $\widehat{\Sigma}$ by
\begin{equation}
\widehat{\Sigma} = \int_0^1 \widehat{G}_Y(t) \lambda_0 (t) dt,
\label{eqn:e659}
\end{equation}
where
\begin{align*}
\widehat{G}_Y(t) & = \frac{1}{n} \sum_{i=1}^n
 Y_i(t) \{ \bwi - \overline{\bm{W}}_Y(t) \}^{\otimes 2},\\
\overline{\bm{W}}_Y(t)  & = 
\frac{n^{-1}\sum_{i=1}^n  Y_i(t)\bwi }{n^{-1}\sum_{i=1}^n  Y_i(t)}.
\end{align*}
Then by just following the arguments on pp.1161-1162 of \cite{HSYYZ2013}
with a sufficiently small $M$, we obtain
\begin{equation}
\widetilde{\Sigma}- \exp \{ - C_XC_g\} \{ 1+ O(L^{-2}) \}
\widehat{\Sigma}\ {\rm is\ n.n.d.}
\label{eqn:e661}
\end{equation}

Finally we recall the definitions of $\overline{\Sigma} $, $\overline{G}_Y(t) $,
and $\bm{\mu}_Y(t) $ in Proposition \ref{prop:prop33}
and set
\begin{equation}
\Delta_3 = \widehat{\Sigma} - \overline{\Sigma}
= - \int_0^1 \overline{Y}(t)
\{ \overline{\bm{W}}_Y(t) - \bm{\mu}_Y(t) \}^{\otimes 2}
\lambda_0 (t) dt
\label{eqn:e665}
\end{equation}
and $\Delta_4 = \overline{\Sigma} - \RE \{ \overline{\Sigma} \}$.
Then we evaluate 
\[
\max_{j,k} | ( \Delta_3 )_{jk} | \quad{\rm and}\quad
\max_{j,k} | ( \Delta_4 )_{jk} |.
\]

Now we give the details for $\Delta_1$, $\Delta_3$, and
$\Delta_4$.

%%%%%%%%%%%%%%%%%%%%%%%%%%%% Delta_1
\noindent
${\bf \Delta_1:}$ We denote the $(jL+l,kL+m)$ element
of $V_n(t,\bgg^*)$ by $v_{jL+l,kL+m}(t)$. Then we have
\begin{equation}
v_{jL+l,kL+m}(t)= (S^{(2)}(t,\bgg^*))_{jL+l,kL+m}
- \frac{(S^{(1)}(t,\bgg^*))_{jL+l}
(S^{(1)}(t,\bgg^*))_{kL+m}}{S^{(0)}(t,\bgg^*)}
\label{eqn:e667}
\end{equation}
and it is easy to see that $|v_{jL+l,kL+m}(t)|$
is uniformly bounded in $j$, $k$, $l$, $m$, and $t$. Besides,
\begin{align}
(S^{(2)}(t,\bgg^*))_{jL+l,kL+m}
& \le C_2 
\begin{cases}
L^{-1},& l=m=1\\
L^{-1/2}|b_l(t)|, & l\ge 2,\ m=1\\
L^{-1/2}|b_m(t)|, & l=1,\ m\ge 2\\ |b_l(t)||b_m(t)|, &
l\ge 2,\ m\ge 2
\end{cases}
\label{eqn:e669}\\
\intertext{and}
(S^{(1)}(t,\bgg^*))_{jL+l}& \le C_3
\begin{cases}
L^{-1/2},& l=1\\ |b_l(t)|, & l\ge 2
\end{cases}.
\label{eqn:e671}
\end{align}

By (\ref{eqn:e667})-(\ref{eqn:e671}) and some calculation, we evaluate
the predictable variation process of $\Delta_1$ and obtain
\begin{equation}
\int_0^1| v_{jL+l,kL+m}(t)|^2 d <\overline{M},\overline{M}>(t)
\le \frac{C_4}{n}
\int_0^1| v_{jL+l,kL+m}(t)|\lambda_0(t) dt \le \frac{C_5}{nL},
\label{eqn:e673}
\end{equation}
where $<\overline{M},\overline{M}>(t) $ is
the predictable variation process of $\overline{M}(t)$.
We used (\ref{eqn:e690}) here.

Thus we have from the exponential inequality for martingales
that
\begin{equation}
\RP \Big(  \max_{j,k}| (\Delta_1)_{jk} | \ge
\frac{x(\log n)^{1/2}}{\sqrt{nL}}\Big)
\le 2(pL)^2\exp \Big\{ -\frac{C_6x^2\log n}{x(\log n)^{1/2} (n^{-1}L)^{1/2} +1}
\Big\}.
\label{eqn:e675}
\end{equation}

%%%%%%%%%%%%%%%%%%%%%%%%%%%% Delta_3
\noindent
${\bf \Delta_3:}$ Notice that $\overline{\Sigma}-\widehat{\Sigma}$ is n.n.d.
Therefore instead of $\Delta_3$, we treat
\begin{align*}
\Delta_3' & = \frac{1}{C_Y}
\int_0^1 \{ \overline{Y}(t) \}^2
\{ \overline{\bm{W}}_Y(t) - \bm{\mu}_Y(t) \}^{\otimes 2}
\lambda_0 (t) dt\\
& = \frac{1}{C_Y}\int_0^1 \Big[ n^{-1}
\sum_{i=1}^n \Big\{ \bwi - Y_i(t)\bm{\mu}_Y(t)
\Big\} \Big]^{\otimes 2}\lambda_0 (t) dt .
\end{align*}

We evaluate $(\Delta_3')_{kl} = (C_Yn^2)^{-1}\sum_{i,j}f_{ij}$,
where $\bm{\mu}_Y(t) =( \mu_{Y1}(t), \ldots, \mu_{Yp}(t))^T$ and 
\[
f_{ij}=\int_0^1 \{ W_{ik}(t) -Y_i(t)\mu_{Yk} (t) \}
\{ W_{jl}(t) -Y_j(t)\mu_{Yl} (t) \}\lambda_0 (t) dt .
\]
Note that $|f_{ij}| \le C_7 L^{-1} $. Thus by applying 
Lemma 4.2 in \cite{HSYYZ2013}, we obtain
\begin{equation}
\RP \Big( \max_{k,l}|(\Delta_3')_{kl}| \ge \frac{x(\log n)^{1/2}}{\sqrt{n}L}
\Big)\le 5(pL)^2\exp 
\Big\{ -\frac{C_8x(n\log n)^{1/2}}{x^{1/2}(n^{-1}\log n )^{1/4}+1} \Big\}.
\label{eqn:e677}
\end{equation}

%%%%%%%%%%%%%%%%%%%%%%%%%%%% Delta_4
\noindent
${\bf \Delta_4:}$ Note that
\[
(\overline{\Sigma})_{kl}= \frac{1}{n}\sum_{i=1}^n
\int_0^1 Y_i(t)\{
W_{ik}(t) -\mu_{Yk} (t) \}
\{ W_{il}(t) - \mu_{Yl} (t) \}\lambda_0 (t) dt
\quad{\rm and}
\]
\[ \Big|\int_0^1 Y_i(t)\{
W_{ik}(t) -\mu_{Yk} (t) \}
\{ W_{il}(t) - \mu_{Yl} (t) \}\lambda_0 (t) dt
% (\overline{\Sigma})_{kl}
\Big|
\le C_9L^{-1}.\]

Applying Bernstein's inequality to $(\overline{\Sigma})_{kl}$, we have
\[
\RP \Big( | ( \Delta_4)_{kl} | \ge \frac{x(\log n)^{1/2}}{\sqrt{n} L}
\Big)\le 2 \exp \Big\{ -\frac{C_{10}x^2\log n}{x(n^{-1}\log n)^{1/2} +1}
\Big\}.
\]
Consequently we have
\begin{equation}
\RP \Big( \max_{k,l}| ( \Delta_4)_{kl} | \ge \frac{x(\log n)^{1/2}}{\sqrt{n} L}
\Big)\le 2 (pL)^2
\exp \Big\{ -\frac{C_{10}x^2\log n}{x(n^{-1}\log n)^{1/2} +1}
\Big\}.
\label{eqn:e679}
\end{equation}

By combining (\ref{eqn:e656}), (\ref{eqn:e657}), (\ref{eqn:e661}), (\ref{eqn:e665})
and (\ref{eqn:e675})-(\ref{eqn:e679}) and exploiting
Lemma \ref{lem:lem73}, we obtain the desired results. 
Hence the proof of the proposition is complete.
\end{pop3}

%%%%%%%%%%%%%%%%%%%%%%%%
%%%%%%%%%%%%%%%%%%%%%%%% Concluding remarks
\section{Concluding remarks}
\label{sec:conclusion}
\noindent
%{\bf To be written later}
%1. Orthogonal decomposition approach for three important models
%2. Numerical studies and promising results
%3. More numerical examination and selection rule for $\lambda$

We proposed an orthonormal basis approach for simultaneous
variable selection and structure identification for varying coefficient
Cox models. We have derived an oracle inequality for the group Lasso procedure
and our method and theory also apply to additive Cox models.
These models are among important structured nonparametric regression models.
This orthonormal basis approach can be used for the adaptive group Lasso
procedure. We presented some preliminary simulation results in this paper.
Extensive numerical examinations and screening-consistent selection rule 
for $\lambda$ are topics of future research.

%%%%%%%%%%%%%%%%%%%%%%%%%%%%%%%%%%%%%%%%%%%%%%%%%%%%%%%%%%% Acknowledgements
\section*{Acknowledgments}
This research is supported by JSPS KAKENHI Grant Number JP 16K05268.

\appendix

\section{Construction and properties of basis functions}
\label{sec:basis}
We describe how to construct $\obt$, 
the properties of $\obt$,
and the approximations to $\bgt$. Set
\[
\Omega_0 = \int_0^1 \bm{B}_0(t) \bm{B}_0^T(t) dt
\quad {\rm and}\quad \overline{\Omega}
= \int_0^1 \obt \overline{\bm{B}}^T(t) dt.
\]

First we describe how to construct $A_0$ and $\obt$. Set
\[
b_1(t) = 1/\sqrt{L} \quad {\rm and} \quad
b_2(t) = \sqrt{12L^{-1}}(t-1/2)\]
and define a inner product on the $L_2$ function space
on $[0,1]$ by
\[
(g_1, g_2) = \int_0^1 g_1(t)g_2(t) dt.
\]
Then we have
\[
\| b_1\|^2 = \| b_2 \|^2 =L^{-1}\quad {\rm and}
\quad (b_1, b_2)=0.
\]
Note that there is some $L$-dimensional vector $\bm{a}_{02}$ satisfying
$ b_2(t) = \bm{a}_{02}^T \bm{B}_0(t)$.

We can obtain $b_j$, $j=3, \ldots, L$, by just applying
the Gram-Schmidt orthonormalization to $(L-2)$ elements of $\bm{B}_0(t)$
with the normalization of $\| b_j\|^2 = L^{-1}$.
Since every $b_j(t)$ is a linear combination of $\bm{B}_0(t)$, we have
\[
\obt = A_0 \bm{B}(t).
\]

Hence we have
\begin{equation}
\overline{\Omega} = A_0 \Omega_0 A_0^T=
\begin{pmatrix}
1/L & \bm{0}^T\\
\bm{0} & \int \bm{B}(t)\bm{B}^T(t)dt
\end{pmatrix}
= \begin{pmatrix}
1/L & \bm{0}^T\\
\bm{0} & A_{-1}\Omega_0A_{-1}^T
\end{pmatrix}
= \frac{1}{L}I.
\label{eqn:e901}
\end{equation}

It is known that for some positive constants $C_1$ and $C_2$,
% $C_3$, and $C_4$,
we have
%\begin{align}
\begin{equation}
\frac{C_1}{L} \le \lambda_{\min}( \Omega_0 ) \le \lambda_{\max}
( \Omega_0 ) \le \frac{C_2}{L} 
\label{eqn:e903}
\end{equation}
%\\ \intertext{and}
%\frac{C_3}{L} & \le \lambda_{\min}( \overline{\Omega} ) \le \lambda_{\max}
%( \overline{\Omega} ) \le \frac{C_4}{L} 
%\label{eqn:e905}
%\end{align}
See Huang et al.\cite{HWZ2004} for more details.
%and Xue and Yang\cite{XY2006} for more details.

Thus (\ref{eqn:e901}) and (\ref{eqn:e903}) imply that
\begin{align}
C_3 & \le \lambda_{\min}( A_0A_0^T) \le \lambda_{\max}
( A_0A_0^T ) \le C_4 
\label{eqn:e906}\\
\intertext{and}
C_5& \le \lambda_{\min}( A_{-1}A_{-1}^T ) \le \lambda_{\max}
( A_{-1}A_{-1}^T ) \le C_6
\label{eqn:e904}
\end{align}
for some positive constants $C_3$, $C_4$, $C_5$, and $C_6$. Note
that (\ref{eqn:e906}) implies that
\[
C_3 \le \lambda_{\min}( A_0^TA_0) \le \lambda_{\max}
( A_0^TA_0 ) \le C_4 .
\]

On the other hand, the definition of $\bm{B}_0(t)$,
(\ref{eqn:e901}), and  (\ref{eqn:e904}) imply
that
\begin{equation}
\int_0^1 b_j (t) dt =0, \ {\rm for}\ j=2,\ldots,L, \quad{\rm and}\quad
\sup_{2\le j \le L} \| b_j \|_\infty = O(1).
\label{eqn:e909}
\end{equation}
Besides, we have for $\bgg_j =( \gamma_{1j}, \bgg_{-1j}^T )^T \in R^L$,
\begin{align}
\bgg_j^T \obt & = \bgg_j^T A_0 \bm{B}_0(t) \quad {\rm and}\nonumber
\\ | \bgg_j^T \obt | & \le
( \bgg_j^T A_0 A_0^T \bgg_j )^{1/2} | \bm{B}_0(t) | \le C_{7} | \bgg_j |
\label{eqn:e913}
%\\ & \le C_{11} | \bgg_j |
\end{align}
uniformly on $[0,1]$ for some positive constant $C_{7}$. Note
that we used (\ref{eqn:e906}) and the local
property of $ \bm{B}_0(t) $ to derive (\ref{eqn:e913}).

Next we consider the approximations to $\bgt $. From Corollary 6.26 in
\cite{Schumaker2007} and Assumption G,
there exist $\bgg_{0j}^* \in R^L$, $j=1,\ldots, p$, satisfying
\begin{equation}
\sum_{j=1}^p \|  g_j - \bm{B}_0^T \bgg_{0j}^* \|_\infty \le
\frac{C_{approx}}{2L^2},
\label{eqn:e914}
\end{equation}
where $C_{approx}$ depends on $C_g$.

In this paper, we use $\obt$ instead of $ \bm{B}_0(t) $. Then
\begin{align*}
\bm{B}_0^T(t) \bgg_{0j}^* & = \overline{\bm{B}}^T(t)
(A_0^T)^{-1}\bgg_{0j}^* = \overline{\bm{B}}^T(t)\overline{\bgg}_{j}^*
\quad ({\rm say})\\
& = \overline{\bm{B}}^T(t) 
\begin{pmatrix}
\overline{\gamma}_{1j}^*\\
\overline{\bgg}_{-1j}^*
\end{pmatrix}\quad ({\rm say}).
\end{align*}

Noticing
\begin{eqnarray*}
\lefteqn{
\sum_{j=1}^p\Big| \int_0^1g_j(t)dt - \frac{\overline{\gamma}_{1j}^*}{L^{1/2}}
- \int_0^1 \overline{\bgg}_{-1j}^{*T}\bm{B}(t) dt \Big|}\\
& = &\sum_{j=1}^p |g_{cj} -L^{-1/2}\overline{\gamma}_{1j}^*|
\le \frac{C_{approx}}{2L^2},
\end{eqnarray*}
we take $\bgg_j^* = 0$ for $\overline{\cal S}_c$,
\begin{align}
\gamma_{1j}^* & = L^{1/2}g_{cj}
\quad {\rm and} \quad \bgg_{-1j}^* = 0
\quad {\rm for}\ j\in {\cal S}_c\cap \overline{\cal S}_n,
\label{eqn:e915}\\
\gamma_{1j}^* & = L^{1/2}g_{cj}
\quad {\rm and} \quad  \bgg_{-1j}^*=
\overline{\bgg}_{-1j}^{*}
\quad {\rm for}\ j\in {\cal S}_n .
\nonumber
\end{align}

Then from (\ref{eqn:e914}),  we have
\begin{align}
& \sum_{j=1}^p \|  g_j - \overline{\bm{B}}^T
\bgg_{j}^* \|_\infty %&
 \le \frac{C_{approx}}{L^2}
\label{eqn:e916}\\
\intertext{and uniformly in $j$,}
\| g_j \|^2 & = |g_{cj}|^2 + \| g_{nj} \|^2 =
\bgg_j^{*T} \overline{\Omega} \bgg_j^* + O(L^{-4})
\nonumber\\
& = \frac{| \gamma_{1j}^* |^2}{L}
+  \bgg_{-1j}^{*T}
\int_0^1 \bm{B}(t) \bm{B}^T(t)dt \bgg_{-1j}^{*}
+ O(L^{-4})\nonumber \\
& =  \frac{| \gamma_{1j}^* |^2}{L}+ \frac{ | \bgg_{-1j}^* |^2}{L}
+O(L^{-4}).\nonumber
\end{align}
We also have
\begin{align}|g_{cj}|^2= \frac{| \gamma_{1j}^* |^2}{L}
\quad {\rm and}\quad
\|g_{nj}\|^2 = \frac{| \bgg_{-1j}^* |^2}{L}
+O(L^{-4}).
\label{eqn:e917}
\end{align}

%\begin{align}\label{eqn:e9}\end{align}
%\smallskip
%\medskip
%\begin{equation}\label{eqn:e9}\end{equation}

\section{Subgradients}
\label{sec:subgradients}
We give $\nabla_j P_1(\bgg) $ and $\nabla_j P_h(\bgg) $ just for reference.

For $\nabla_j P_1(\bgg) $, we have
\begin{align*}
\nabla_j |\gamma_{1j}| & =
\begin{cases} \sign ( \gamma_{1j} ), & |\gamma_{1j}| \ne 0\\
\epsilon_{1j}, & \gamma_{1j}=0
\end{cases}\\
\intertext{and}
\nabla_j |\bgg_{-1j}| & =
\begin{cases} \bgg_{-1j}/|\bgg_{-1j}|, & |\bgg_{-1j}| \ne 0\\
\bm{\epsilon}_{-1j}, & \bgg_{-1j}=0
\end{cases},
\end{align*}
where $|\epsilon_{1j}| \le 1$ and $ |\bm{\epsilon}_{-1j}| \le 1$.

Next we deal with $\nabla_j P_h ( \bgg) $. Recall that
\[
\nabla_j P_h(\bgg) = \nabla_j ( |\gamma_{1j}|^q + |\bgg_{-1j}|^q )^{1/q}
+  \nabla_j |\bgg_{-1j}| .
\]

Set
\[
\nabla_j ( |\gamma_{1j}|^q + |\bgg_{-1j}|^q )^{1/q} =
\begin{pmatrix}
d_{1j}\\
\bm{d}_{-1j}
\end{pmatrix},
\]
where $d_{1j} \in R$ and $ \bm{d}_{-1j} \in R^{L-1}$.

When $|\gamma_{1j}|=0$ and $ |\bgg_{-1j}| =0$,
\[
d_{1j} = \epsilon_{1j}\quad {\rm and}\quad
\bm{d}_{-1j}= \bm{\epsilon}_{-1j},
\]
where $|\epsilon_{1j} |\le a$ and
$|\bm{\epsilon}_{-1j} |\le b$ such that
$(a,b)$ satisfies $(1+t^q)^{1/q} \ge a+bt$ for any $t\ge 0$.
This follows from the definition of subgradient and we note
that $0\le a \le 1$ and $0\le b \le 1$.

When $|\gamma_{1j}|\ne 0$ and $ |\bgg_{-1j}| =0$,
\[
d_{1j} = \sign ( \gamma_{1j} ) \quad {\rm and}\quad
\bm{d}_{-1j}= 0.
\]

When $|\gamma_{1j}|=0$ and $ |\bgg_{-1j}| \ne 0$,
\begin{equation}
d_{1j} = 0 \quad {\rm and}\quad \bm{d}_{-1j}=
\bgg_{-1j}/|\bgg_{-1j}|.
\label{eqn:e919}
\end{equation}
This property is essential to  hierarchical selection
for $g_{cj}$ and $g_{nj}(t)$. See \cite{ZRY2009}.

When $|\gamma_{1j}| \ne 0$ and $ |\bgg_{-1j}| \ne 0$,
\begin{align*}
d_{1j} & =( |\gamma_{1j}|^q + |\bgg_{-1j}|^q )^{\frac{1}{q}-1}
\sign ( \gamma_{1j} )| \gamma_{1j}|^{q-1}
\\
\intertext{and}
\bm{d}_{-1j} & = ( |\gamma_{1j}|^q + |\bgg_{-1j}|^q )^{\frac{1}{q}-1}
\frac{\bgg_{-1j}}{|\bgg_{-1j}|}|\bgg_{-1j}|^{q-1}.
\end{align*}

%\begin{align}\label{eqn:e9}\end{align}
%\smallskip
%\medskip
%\begin{equation}\label{eqn:e9}\end{equation}

\section{Proofs of technical lemmas}
\label{sec:technical}
\begin{pol1} From the definitions of $\obt$ and $\bwi$.
we have
\begin{equation}
\bt^T ( \bwi - \bm{W}_j(t) )
= \bt^T ( I_p \otimes A_0) ( \bxi \otimes \bm{B}_0(t)
-\bm{X}_j(t) \otimes \bm{B}_0(t) ).
\label{eqn:e921}
\end{equation}
Notice that for $\bt = (\bt_1^T, \ldots, \bt_p^T)^T$,
\begin{equation} | \bt_k^T A_0
\bm{B}_0(t)| \le  |  A_0^T\bt_k | \le \{ \lambda_{\max} 
( A_0A_0^T ) \}^{1/2} |\bt_k |.
\label{eqn:e923}
\end{equation}
Here we used that $|\bm{B}_0 (t) |\le 1$.

Consequently (\ref{eqn:e921}) and (\ref{eqn:e923}) yield that
\begin{eqnarray*}
\lefteqn{ | \bt^T ( \bwi - \bm{W}_j(t) ) | }\\
& \le & \sum_{k=1}^p |X_{ik}(t) - X_{jk}(t)| | \bt_k^T A_0
\bm{B}_0(t)| \\
& \le & 2C_X\{ \lambda_{\max}(  A_0A_0^T) \}^{1/2}
\sum_{k=1}^p | \bt_k | \, \le \, C_W P_1 (\bt ).
\end{eqnarray*}

Hence the proof is complete.
\end{pol1}

\begin{pol2} This lemma is just a version of Lemma 3.2 in \cite{HSYYZ2013}.
We can verify this lemma in the same way by taking
\[
a_i(t) =\bt^T \{ \bwi - \widetilde{\bm{W}}_n(t, \bgg^* )\}
\quad {\rm and}\quad w_i(t) = Y_i(t) \exp \{ \bgg^{*T}\bwi \}
\]
in the proof. The details are omitted.
Hence the proof is complete.
\end{pol2}

\begin{pol3} This is almost proved in \cite{HSYYZ2013}.
We should just note that
\begin{align*} | \bgg^T (\Sigma_1 -\Sigma_2) \bgg |
& \le |\bgg|_1^2 \max_{j,k} | ( \Sigma_1 -\Sigma_2 )_{jk} |
\le L \{ P_1(\bgg ) \}^2 \max_{j,k} | ( \Sigma_1 -\Sigma_2 )_{jk} |,
\\
P_1 (\bgg ) & \le (1+\zeta ) P_1 (\bgg_{\cal S} ),\quad {\rm and}\quad
P_1 (\bgg_{\cal S} ) \le s_0^{1/2} | \bgg |.
\end{align*}

When $\Sigma_2 - \Sigma_1 $ is n.n.d., we have
\[ | \bgg^T (\Sigma_1 -\Sigma_2) \bgg |
\le \bgg^T \Delta \bgg \le L \{ P_1(\bgg ) \}^2
\max_{j,k} | ( \Delta )_{jk} |.
\]

Hence the proof is complete.
\end{pol3}
%\begin{align}\label{eqn:e9}\end{align}
%\smallskip
%\medskip
%\begin{equation}\label{eqn:e9}\end{equation}

\section{Derivatives of the B-spline basis}
\label{sec:bspline}
In this section, we examine properties of
\[
\int_0^1 \bm{B}_0'(t) (\bm{B}_0'(t))^T dt
\]
and describe why we have adopted the orthogonal decomposition approach while
the other authors have considered the $L_2$ norm of the estimated derivatives
when they deal with structure identification for additive models
or partially linear additive models.

We take a function $g_A(t)$ on $[0,1]$ defined by
\[
g_A(t) = \sin (2\pi A t)
\]
for $A \to \infty$ sufficiently slowly. Then it is easy to see 
\[
\| g_A \|^2 \sim 1,
\quad\| g_A' \|^2 \sim A^2,
\quad{\rm and}\quad  
\quad\| g_A'' \|^2 \sim A^4.
\]
On the other hand, we can approximate this $g(t)$ by $\bm{B}_0(t)
\bm{\gamma}_A$ accurately enough and we have
\[
\bm{\gamma}_A^T \Omega_0 \bm{\gamma}_A \sim 1,
\quad |\bm{\gamma}_A |^2 \sim L,
\quad {\rm and} \quad
\bm{\gamma}_A^T \int_0^1 \bm{B}_0'(t) (\bm{B}_0'(t))^T dt
\bm{\gamma}_A \sim A^2\to \infty.
\]

This means some eigenvalues of
$\int_0^1 \bm{B}_0'(t) (\bm{B}_0'(t))^T dt$
have the order larger than $L^{-1}$. Hence we cannot follow the proofs
in the papers based on the $L_2$ norm of the estimated derivatives
because the present eigenvalue property violates their assumptions on
matrices similar to
\[
\int_0^1 \bm{B}_0''(t) (\bm{B}_0''(t))^T dt.
\]
The above matrix also should have some larger eigenvalues. Besides,
it is more difficult to estimate the derivatives of the coefficient
functions. This is why we have adopted the orthogonal decomposition
approach. Zhang et al.\cite{ZCL2012} is based on the smoothing spline method
and it is difficult to apply their ingenious approach to the loss function
other than the $L_2$ loss function.

%\begin{align}\label{eqn:e9}\end{align}
%\smallskip
%\medskip
%\begin{equation}\label{eqn:e9}\end{equation}

\section{Proofs for other models}
\label{sec:additionalproofs}
We outline necessary changes in the proofs for the former model
in section \ref{sec:other} since both models in the section can
be treated in almost the same way as the time-varying coefficient
model. Especially, almost no change is necessary to the proofs
of Proposition \ref{prop:prop31} and Theorem \ref{thm:thm31}.

We assume standard assumptions for varying coefficient models here.

\medskip
\noindent
Proof of Proposition \ref{prop:prop32}) The poof consists
of (\ref{eqn:e646})-(\ref{eqn:e649}) and (\ref{eqn:e651})-(\ref{eqn:e658}).

\smallskip
\noindent
(\ref{eqn:e646})-(\ref{eqn:e649}): Note that $|b_{0j}(t)|$ is replaced with
$n^{-1}\sum_{i=1}^n|b_{0j}(Z_{i}(t))|$. When we evaluate the predicable
variation process in (\ref{eqn:e645}), 
\[
\int_0^1 |b_{0j}(t)|^2 \lambda_0(t) dt \le C \int_0^1 |b_{0j}(t)| \lambda_0(t) dt
\]
is replaced with
\begin{equation}
\int_0^1 \Big\{ n^{-1}\sum_{i=1}^n|b_{0j}(Z_{i}(t))| \Big\}^2  \lambda_0(t) dt
\le C \int_0^1
n^{-1}\sum_{i=1}^n b_{0j}^2(Z_{i}(t)) \lambda_0(t) dt.
\label{eqn:ee001}
\end{equation}

We can evaluate the second term in (\ref{eqn:ee001}) by using Bernstein's inequality
and
\[
\RE \Big\{ n^{-1}\int_0^1
\sum_{i=1}^n b_{0j}^2(Z_{i}(t)) \lambda_0(t) dt  \Big\}
= \int_0^1 \RE \{ b_{0j}^2(Z_1(t)) \} \lambda_0(t) dt
= O(L^{-1}).
\]

\smallskip
\noindent
(\ref{eqn:e651})-(\ref{eqn:e658}): When we apply the martingale
exponential inequality, (\ref{eqn:e690}) is replaced with
\[
\frac{1}{n} \sum_{i=1}^n \int_0^1 b_j^2(Z_i(t))\lambda_0 (t) dt.
\]
We can evaluate this expression by using Bernstein's inequality
and
\begin{align*}
\RE \Big\{ \int_0^1 b_j^2(Z_1(t))\lambda_0 (t) dt \Big\}
& \le C\bm{a}_{0j}^T \int_0^1 \RE \{ \bm{B}_0(Z_1(t))(
\bm{B}_0(Z_1(t)))^T \}
\lambda_0 (t) dt\bm{a}_{0j} \\
& = O(L^{-1}).
\end{align*}
We need some assumptions for $\RE \{ \bm{B}_0(Z_1(t))(
\bm{B}_0(Z_1(t)))^T \}$ as for $\Omega_0$ in \ref{sec:basis}.

\medskip
\noindent
Proof of Proposition \ref{prop:prop33}) The proof
consists of evaluating $\Delta_1$, $\Delta_3$,
and $\Delta_4$.

\smallskip
\noindent
${\bf \Delta_1}$: We should just follow the line of (\ref{eqn:e651})-(\ref{eqn:e658}).

\smallskip
\noindent
${\bf \Delta_3}$: This is almost a U-statistic and we can also apply 
the exponential inequality for U-statistics as (3.5) in \cite{GLZ2000}
to the part of a U-statistic.

\smallskip
\noindent
${\bf \Delta_4}$: This is a sum of bounded independent random variables
and we can deal with this by applying Bernstein's inequality.

%\begin{align}\label{eqn:e9}\end{align}
%\smallskip
%\medskip
%\begin{equation}\label{eqn:e9}\end{equation}

%% References


%\bibliographystyle{elsarticle-num}



%%%%%%%%%%%%%%%%%%%%%%%%%%%%%%%%%%%%%%%%%%%%%%%%%%%%%%%%%%%%%%%%%%%%%% References

\end{document}